\documentclass[preprint,12pt]{elsarticle}
\usepackage{overpic}
\usepackage{bm}
\usepackage{xcolor}
\usepackage{amsmath}
\usepackage{cases}
\usepackage{float}
\usepackage{amsthm}
\usepackage{geometry}
\usepackage{graphicx}
\usepackage{tikz}
\usepackage{algorithm}
\usepackage{algpseudocode}
\usepackage{amsfonts,amssymb}
\usepackage{caption}
\usepackage{tikz}
\usepackage{algorithm}
\usepackage{algorithmicx}
\usepackage{algpseudocode}
\usepackage{threeparttable}
\biboptions{sort&compress}
\usepackage{subfig}
\usepackage{booktabs}
\usepackage{array, caption, threeparttable}
\usepackage[font=small,labelfont=bf,labelsep=none]{caption}

\begin{document}
	\begin{frontmatter}
		\title{Mass-preserving spatio-temporal adaptive PINN for Cahn-Hilliard equations with strong nonlinearities and singularities}
		
		\author[1]{Qiumei Huang}
		\ead{qmhuang@bjut.edu.cn}
		
		\author[3]{Jiaxuan Ma}
		\ead{majiaxuan@emails.bjut.edu.cn}
		
		\author[1]{Zhen Xu\cormark[2]}
		\ead{xuzhenmath@bjut.edu.cn}
		\address[1]{School of Mathematics, Statistics and Mechanics, Beijing University of Technology, Beijing 100124, China}
        \address[3]{Xicheng District Bureau of Statistics, Beijing 100032, China}
	\cortext[2]{Corresponding author.}	

		\begin{abstract}
   			As one kind of important phase field equations, Cahn-Hilliard equations involve high-order spatial derivatives, strong nonlinearities, and even solution singularities when certain bulk potentials (e.g., Flory-Huggins) are used. When using the physics-informed neural network (PINN) to simulate the long time evolution of the solution, it is necessary to decompose the time domain to capture the transition of solutions in different time. Moreover, the standard PINN can not maintain the mass conservation property for the equations exactly. We propose a novel mass-preserving spatio-temporal adaptive PINN, which adaptively divides the time domain according to the rate of energy decrease, and solves the Cahn-Hilliard equation within each subinterval. To improve the prediction accuracy, spatial adaptive sampling is employed in the subdomain to select points with large residual value which are added to the training samples. 
			Notably, a mass constraint is added to the loss function to compensate for the mass degradation problem of the PINN method when solving Cahn-Hilliard equations. Numerical experiments are presented to illustrate the effectiveness of the proposed method in solving complex phase-fields models, including the Cahn-Hilliard equations with different bulk potentials, the three-dimensional Cahn-Hilliard equation with singularities, and the system of Cahn-Hilliard equations. 		 
		\end{abstract}
	\end{frontmatter}
	\section{Introduction}\par
	The phase field model is widely applied to a variety of problems in materials science and life science. It utilizes the order parameter $u$ to capture the diffuse interface between two phases. One of the typical and classical equations of the phase field model is the Cahn-Hilliard equation proposed by Cahn and Hilliard \cite{CahnHilliard1958,Cahn1961}, which was originally used to simulate phase separation phenomena in binary metal mixtures. Then it is widely extended to simulate the motion of free interfaces during phase transitions. 
	A variety of efficient numerical methods related to time discretization have been proposed for the phase field equations, such as convex splitting schemes \cite{Eyre1998,Wise2009,Chen2011,Gomez2011,Guillen2013}, stabilized schemes \cite{FengTangYang2013,Tang2007,Shen2010,XuTang2006}, Invariant Energy Quadratization (IEQ) method and its variants \cite{Yang2016,Xu2019,Pan2023}, Scalar Auxiliary Variable (SAV) and its variants \cite{ShenXuYang2018,ShenXu2018,Li2019,Liu2020,Jiang2022}, exponential time differencing (ETD) schemes \cite{Ju2015,Ju2019} and so on. 
    However, these methods encounter certain difficulties when the energy of the phase field model is extremely complicated, where the equation contains singularities or complex boundary conditions. Moreover, the curse of dimensionality remains a challenge for traditional numerical methods in high-dimensional spaces. \par

Deep learning offers a potential panacea to solve the high dimensional problems, so it has attracted the attention of researchers and then has been widely utilized in various fields with the improvement of computational resources and the iteration of algorithms. Physics-informed neural networks (PINNs), proposed by Raissi et al.~\cite{RaissiKarniadakis2019}, embed governing physical laws and encodes the underlying  PDEs directly into the loss function, enabling the solution of forward and inverse problems for PDEs. 
For the phase field equation, PINN offers a mesh-free framework that avoids delicate spatial discretization, handles nonlinearities naturally via automatic differentiation and bounded output mappings, and allows physical constraints such as mass conservation to be incorporated directly into the loss function. These features make PINN a natural candidate for phase-field modeling.\par

Despite their success, standard PINNs face limitations when solving complex problems such as the Cahn-Hilliard equation, which involves high-order spatial derivatives, strong nonlinearities, singularities, and stiffness. 
In particular, a single neural network often fails to capture the fast initial phase separation followed by slow coarsening over long time intervals, and the baseline PINN does not inherently preserve mass conservation. To address these issues, various improvements have been proposed, including adaptive sampling strategies (e.g., RAR \cite{Lu2021}, DAS-PINN \cite{Tang2023}, FI-PINN \cite{GaoZhou2023}, TCAS-PINN\cite{Guo2022}, RAD \cite{Wu2023}), time-domain decomposition (e.g., Zhao and Wight  \cite{Zhao2021}, PPINN \cite{Meng2020}), adaptive loss weighting \cite{Wang2021,Xiang2022}, and sampled-data control strategies for neural networks \cite{Chandrasekar2025,Ma2024}. Additionally, conservation constraints have been incorporated into the loss function to enforce mass conservation \cite{Lin2022}. 

	In this work, we propose the mass-preserving spatio-temporal adaptive PINN for the complex Cahn-Hilliard equations. The fourth-order Cahn-Hilliard equation is converted to a second-order system of equations to reduce computational costs, and the system is then simulated by adaptively partitioning the time domain to capture the rapid transitions of phase patterns over time. The length and number of time subintervals are adjusted according to the rate of energy decrease, and within each subinterval, spatial adaptive sampling is employed to enhance training effectiveness. 
    A soft mass constraint is added to the loss function to preserve mass conservation. Unlike the uniform time used in \cite{Zhao2021}, our adaptive time partitioning based on energy decrease significantly reduces the absolute errors for the Cahn-Hilliard equation with Ginzburg-Landau potential.
    For the Flory-Huggins potential, which exhibits high singularity, we add a mapping to the network output to ensure stability.  
    Finally, we simulate the constrained binary blended polymer model to demonstrate the effectiveness of our method for complex problems.
    The numerical results show that the predicted solutions maintain mass conservation well. In addition, we compare our method with representative PINN variants, including conservative PINN \cite{Jagtap2020cPINN} and PINN with adaptive activation functions \cite{Jagtap2020Adaptive}, to further demonstrate the advantages of the proposed methods. \par

	The rest of this article is organized as follows. In Section 2, we provide a brief overview of the baseline PINN method, followed by the relevant improvement methods and strategies. Finally, we present the mass-preserving spatio-temporal adaptive PINN for high-order phase field equations with mass conservation. In Section 3, we present several numerical examples to demonstrate the accuracy and effectiveness of our algorithm, including simulations of the three-dimensional Cahn-Hilliard equation with Flory-Huggins potential and the binary blended polymers model, as well as additional comparative and sensitivity studies. Finally, we give conclusions in Section 4. \par

	\section{ Methodology}
	In this section, we first briefly review the Cahn-Hilliard equations. Then we introduce the baseline PINN and some preliminary improvements to the PINN, including the sampling method and the loss function settings. Finally, we improve and develop our spatio-temporal adaptive PINN method for high-order phase field equations with mass conservation.

	\subsection{Cahn-Hilliard equations}
	Consider the following energy functional:
	\begin{equation}\label{E}
		E(u)=\int_{\Omega}\bigg(\frac{\varepsilon^2}{2}|\nabla u|^2+F(u)\bigg)d\bm{x},
	\end{equation}
	$u(\bm{x},t)$ is the phase variable, $\bm{x}\in\Omega\subseteq\mathbb{R}^d,t\in[0,T]$, $\varepsilon$ is a positive constant associated
	with the diffuse interface width, $F(u)$ is the nonlinear bulk
	potential. Two common types of nonlinear bulk potentials are shown as follows \cite{CahnHilliard1958,Copetti1992,Elliott1996,Liu2003}:\\
	$\bullet$ Ginzburg-Landau double-well type potential:
	\begin{equation}
		F_{db}(u)=\frac{1}{4}(u^2-1)^2,\quad u\in(-\infty,+\infty).
	\end{equation}
	$\bullet$ Logarithmic Flory-Huggins potential:
	\begin{equation}\label{fh_energy}
		F_{fh}(u)=(1+u)\ln(1+u)+(1-u)\ln(1-u)-\frac{\theta}{2}u^2,\quad u\in(-1,1).
	\end{equation} \par
	The Cahn-Hilliard equations are then obtained under the $H^{-1}$ gradient flow, as follows,
	\begin{equation}\label{CH}
		u_t=\Delta(-\varepsilon^2\Delta u+f(u)),
	\end{equation}
	where $f(u)=F'(u)$. Periodic boundary conditions and Neumann boundary conditions are used on area $\Omega$. As mentioned before, Cahn-Hilliard equation has the properties of mass conservation and energy dissipation:\par
	$\bullet$ mass conservation
	\begin{equation}
		\int_{\Omega}u(\textbf{x},t)d\textbf{x}=\int_{\Omega}u(\textbf{x},0)d\textbf{x},\quad\forall t>0.
	\end{equation}\par
	$\bullet$ energy dissipation
	\begin{equation}
		E(u(t))\leq E(u(s)),\quad \forall t\geq s.
	\end{equation}\par
    The Cahn-Hilliard equation has high nonlinearity and even singularity, and it has spatial higher order derivatives, which make the computation difficult. In addition, it is also noteworthy whether the predicted solutions satisfy the physical properties of the Cahn-Hilliard equation
    \subsection{Baseline PINN and some preliminary improvements}\label{modified_Pinn}
	The baseline PINN can efficiently solve partial differential equations with smooth solutions. However, it is inadequate for problems involving sharp interfaces or long-time integration of time-dependent PDEs. 
    More importantly, it cannot inherently preserve the physical properties such as mass conservation. In this section, we introduce two preliminary improvements to address these limitations: adaptive sampling in space and time (Section 2.2.2) and a mass constraint imposed on the loss function (Section 2.2.3).

	\subsubsection{Baseline PINN}\label{Pinn}
	Consider the following Cahn-Hilliard equation:
	\begin{equation}\label{nonlinerPDE}
		\begin{split}
			&u_t+\varepsilon^2\Delta^2u-\Delta(u^3-u)=0,\quad \bm{x}\in[-1,1],\quad t\in[0,T],\\
			&u(\bm{x},0)=u_0(\bm{x}),\quad \bm{x}\in\Omega,\\
			&u(-1,t)=u(1,t),\quad t\in[0,T].\\
		\end{split}
	\end{equation}
	Here, $\bm{x}$ and $t$ are the space and time coordinates respectively, $u$ is the solution of Equation (\ref{nonlinerPDE}).\par
	\begin{figure}[H]
		\centering
		\begin{tikzpicture}[scale=0.7]
			\tikzstyle{every node}=[font=\small,scale=0.9,align=center]
			\foreach \x in{3,1,-3}
			\fill[green!70](-3,\x)circle(5pt)node(a\x){};
			\node(a-1)at(-3,-1){$\vdots$};
			\foreach \x in{-2,-4,2,4}
			\fill[purple!60](1,\x)circle(5pt)node(b\x){};
			\node(b0)at(1,0){$\vdots$};
			\foreach \x in{-2,-4,2,4}
			\fill[purple!60](5,\x)circle(5pt)node(c\x){};
			\node(c0)at(5,0){$\vdots$};
			\fill[orange!80](9,0)circle(5pt)node(n){};
			
			\node(y1)at(-3.5,1){$x_1$};
			\node(y1)at(-3.5,-3){$x_n$};
			\node(y2)at(-3.5,3){$t$};
			
			\node at(-3,4.5){Input\quad layer};
			\node at(3,4.5){Hidden\quad layer};
			\node at(9,4.5){Output\quad layer};
			\node(o)at(10,0){$\hat{u}_\theta(x,t)$};
			\foreach \x in{3,1,-1,-3}
			{\foreach \y in{-2,-4,0,2,4}
				{\draw(a\x)--(b\y);
				}
			}
			\foreach \x in{-2,-4,0,2,4}
			{\foreach \y in{-2,-4,0,2,4}
				{   
					\draw(b\x)--(c\y);
					\draw(c\y)--(n.west);
					
				}
			}
		\end{tikzpicture}
		\caption{Fully connected neural network}
		\label{neural_network}
	\end{figure}
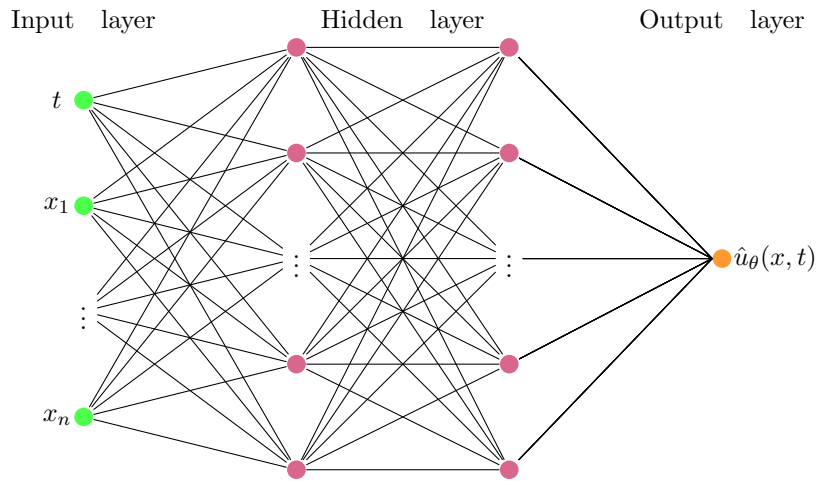
	The main idea of PINN is to use the deep neural network $\hat{u}_{\theta}(\bm{x},t)$ as shown in Fig.~\ref{neural_network} with parameters $\theta$ to approximate the solution $u(\bm{x},t)$ of PDEs. We consider the feedforward fully connected neural network $N_L(\bm{x})$ with $L-1$ hidden layers and $n$ neurons per layer, defined as follows:
	\begin{equation}
		\begin{split}
			\text{input layer}:&\quad N_0(\bm{x})=\bm{x},\\
			\text{hidden layer}:&\quad N_i(\bm{x})=\sigma(\bm{W}_i N_{i-1}(\bm{x})+\bm{b}_i),\quad i=1,\cdots,L-1,\\
			\text{output layer}:&\quad N_L(\bm{x})=\bm{W}_L N_{L-1}(\bm{x})+\bm{b}_L,\\
		\end{split}
	\end{equation}
	where $\sigma$ is the nonlinear activation function. $\bm{W}_i$ and $\bm{b}_i$ represent the weights and biases at $i$-th layer. $\theta=\{\bm{W}_i,\bm{b}_i\}_{i=1}^L$ represents the trainable parameters of the neural network. \par

	In PINN, the gap between predicted solution and real solution is measured by the following loss function:
	\begin{equation}
		\label{Loss_basic}
		\mathcal{L}(\theta)=\lambda_{u}\mathcal{L}_{u}(\theta)+\lambda_{b}\mathcal{L}_{b}(\theta)+\lambda_{r}\mathcal{L}_{r}(\theta).
	\end{equation}
	where,
	\begin{equation}
		\begin{split}
			&\mathcal{L}_{u}(\theta)=\frac{1}{N_{u}}\sum_{i=1}^{N_{u}}|\hat{u}_\theta(\bm{x}_{u}^i,0)-u_0(\bm{x}_{u}^i)|^2,\\
			&\mathcal{L}_{b}(\theta)=\frac{1}{N_{b}}\sum_{i=1}^{N_{b}}|\hat{u}_\theta(-1,t_{b}^i)-\hat{u}_\theta(1,t_{b}^i)|^2,\\
			&\mathcal{L}_{r}(\theta)=\frac{1}{N_{r}}\sum_{i=1}^{N_{r}}\bigg|\frac{\partial\hat{u}_\theta}{\partial t}(\bm{x}_{r}^i,t_{r}^i)+\varepsilon^2\Delta^2\hat{u}_\theta(\bm{x}_{r}^i,t_{r}^i)-\Delta(\hat{u}^3_\theta(\bm{x}_{r}^i,t_{r}^i)-\hat{u}_\theta(\bm{x}_{r}^i,t_{r}^i))\bigg|^2.
		\end{split}
	\end{equation}
	Here, $\{\bm{x}_u^i,0\}_{i=1}^{N_u}$ denotes the initial points and $\{-1,t_b^i\}_{i=1}^{N_b}$, $\{1,t_b^i\}_{i=1}^{N_b}$ denote the boundary points. $\{\bm{x}_r^i,t_r^i\}_{i=1}^{N_r}$ denotes the set of residual points sampled from $[-1,1]\times[0,T]$ by random sampling methods like Latin hypercube sampling (LHS). The derivatives in the loss function are calculated through automatic differentiation in neural networks \cite{RaissiKarniadakis2019}, which avoids complex discrete derivative approximations. The weights $\lambda_u$, $\lambda_b$ and $\lambda_r$ are penalty factors to balance the various loss terms and can be adjusted to obtain the best training result.\par
	PINN transforms the problem of solving PDEs into a neural network optimization problem. The goal of training PINN is to find the optimal parameters $\theta^*$ that minimize the loss function, and this process is implemented using optimization algorithms such as stochastic gradient descent \cite{Theodoridis2015}, Adam \cite{King2014} and L-BFGS \cite{Liu1989}. The numerical experiments in this paper are conducted by initially training with the Adam optimizer, followed by fine-tuning using the L-BFGS optimizer.
    
	\subsubsection{Adaptive Sampling in Space and Time}\label{Sampling}
PINN can be regarded as a kind of semi-supervised learning. Simulation of the initial and boundary conditions belongs to supervised learning. The unsupervised learning part is the PDE residual:
\begin{equation}\label{loss_f}
\mathcal{L}_{r}(\theta)=\frac{1}{N_{r}}\sum_{i=1}^{N_{r}}\bigg|\frac{\partial\hat{u}_\theta}{\partial t}(\textbf{x}_{r}^i,t_{r}^i)+N(\hat{u}_{\theta}(\textbf{x}_{r}^i,t_{r}^i))\bigg|^2.
\end{equation}
Since the PDE residual loss is evaluated on discrete residual points, the location and distribution of these residual points are critical to the performance of PINN. In \cite{Wu2023}, researchers summarized three types of strategies for selecting residual points in training: fixed residual points, uniform and nonuniform points with resampling. We introduce our non-uniform adaptive sampling method used in this paper.

\begin{figure}[H]
		\centering
		\includegraphics[scale=0.5]{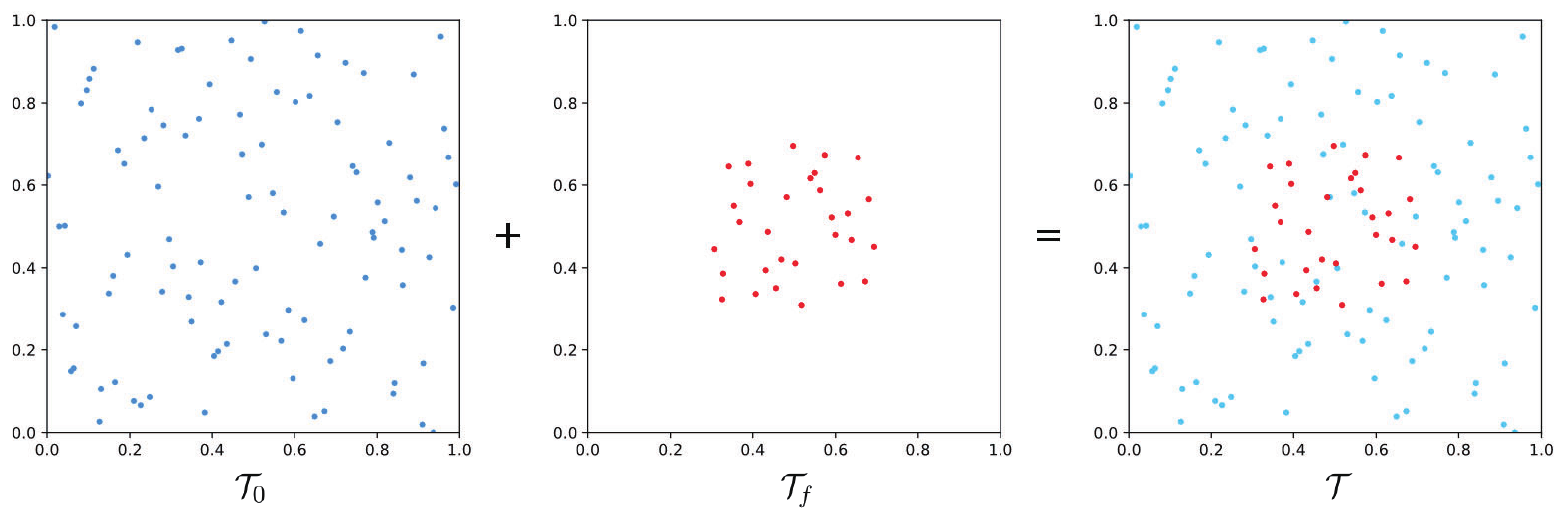}
		\caption{Adaptive sampling in space}
		\label{space_sampling}
	\end{figure}

The PDE residual (\ref{loss_f}) is evaluated at residual points, and then the value of PDE residual can be used as an indicator for selecting residual points. We first periodically stop training to select residual points corresponding to large PDE residuals. Then, these points are added to the set of residual points, shown in Fig.~\ref{space_sampling}. The steps of adaptive sampling in space are illustrated in Algorithm~\ref{alg2}. This method makes more residual points fall into the parts with larger PDE residuals and helps PINN achieve desirable computational results \cite{Lu2021,Zhao2021}.

\begin{figure}[H]
		\centering
		\includegraphics[scale=0.55]{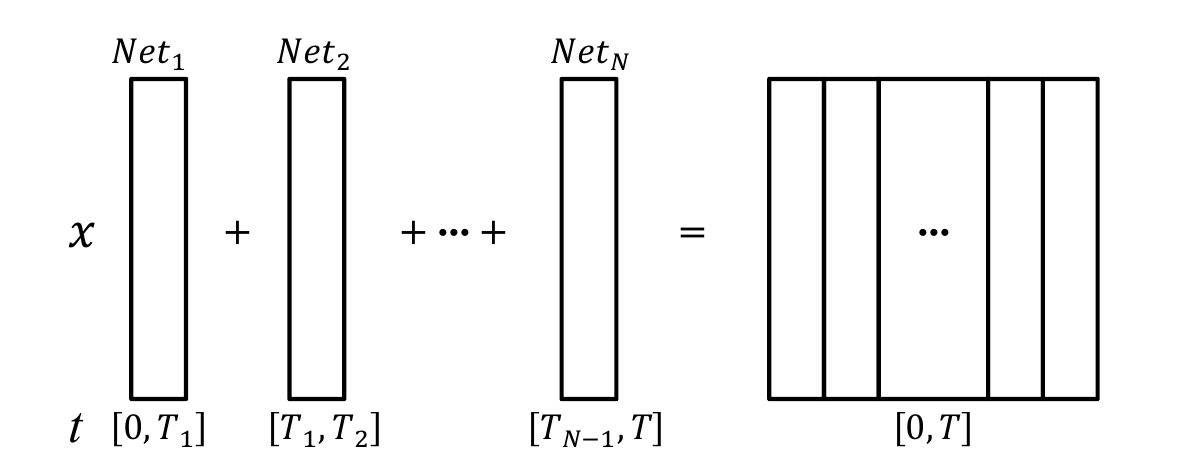}
		\caption{Adaptive sampling in time}
		\label{time_sampling}
	\end{figure}

Considering the sharp transition of the phase pattern with time, we use PINN to solve the Cahn-Hilliard equation sequentially in several time intervals, shown in Fig.~\ref{time_sampling}. The full time domain is generally divided equally. Since the sharp transition of the phase pattern is accompanied by a rapid decrease of energy, we adaptively adjust the length and number of time segments according to the energy decrease rate, which makes the time subintervals unequal.

	\par

	\begin{algorithm}[H]\footnotesize                 
		\caption{Adaptive sampling in space}          
		\label{alg2}      
		\begin{algorithmic}[1] 
			\State Initialize the parameters of PINN with loss function:$\mathcal{L}=\mathcal{L}_u+\mathcal{L}_b+\mathcal{L}_r$;
			\State Sample the initial,boundary and residual points $\mathcal{T}_0$ using LHS, let $\mathcal{T}=\mathcal{T}_0$;
			\State Train the PINN for a certain number of iterations;
			\State Select $m$ points with large residual values as $\mathcal{T}_f$, and let $\mathcal{T}=\mathcal{T}+\mathcal{T}_f$;
			\State Train the neural network with the points set $\mathcal{T}$;
			\State Calculate the mass error;
			\Repeat
			\State Step $4\sim 6$;
			\Until{The total number of iterations or $\mathcal{L}$ is no longer reduced}\\
			\Return $\hat{u}_\theta(x,t)$;
		\end{algorithmic}
	\end{algorithm}
To determine the partition points quantitatively, we monitor the energy gradient at candidate time nodes. Let $|E'(t_k)|$ denote the energy gradient at a candidate node $t_k$, and $|E'(0)|$ the initial gradient. A node is selected as a partition point if
\[
|E'(t_k)| > \alpha \cdot |E'(0)|,
\]
where $\alpha \in (0,1)$ is a threshold ratio that controls the sensitivity of the time partitioning. 
A smaller $\alpha$ leads to a finer partition (more subintervals), while a larger $\alpha$ yields a coarser one. The sensitivity of $\alpha$ is discussed in Section 3.1.

In summary, the proposed spatio-temporal adaptive sampling offers two key advantages. Spatially, residual-based resampling focuses on collocation points near moving interfaces, improving accuracy per training epoch. Temporally, partitioning the time domain according to the energy decrease rate resolves fast initial dynamics while avoiding over-resolving the slow coarsening regime. Unlike existing adaptive PINNs that rely on uniform time subdivision or loss-threshold-based stepping, our time partition is physically driven by the energy dissipation rate - a novel contribution of this work. The combination yields a truly adaptive PINN that significantly enhances both accuracy and efficiency for Cahn-Hilliard equations.

	\subsubsection{Mass Constraint}\label{mass_conserivate}
	
In the experiments to solve the Cahn-Hilliard equations, we found that the mass of the predicted solution by the baseline PINN does not preserve mass conservation, as shown in Fig.~\ref{mass_1}. This is because computing higher-order derivatives by backpropagation is costly for neural networks. We reduce the order of equation \eqref{CH} by introducing a new variable and set the output layer of the neural network to two neurons that output the predicted values of $\varphi$ and $u$, denoted as $\hat{\varphi}_\theta$ and $\hat{u}_\theta$, respectively. After that, the residual loss term of the equation takes the following form:
	\begin{equation}\label{L_r}
		\mathcal{L}_{r}(\theta)=\frac{1}{N_{r}}\sum_{i=1}^{N_{r}}\bigg(\bigg|\frac{\partial\hat{u}_\theta}{\partial t}-\Delta\hat{\varphi}_\theta\bigg|^2+|\hat{\varphi}_\theta-(-\varepsilon^2\Delta \hat{u}_\theta+f(\hat{u}_\theta))|^2\bigg).
	\end{equation}
We note the first term in (\ref{L_r}) can only be optimized to be as small as possible while not exactly equal to $0$ after the neural network training process. This then results in mass degradation for the Cahn-Hilliard equation. \par

To overcome this problem, we add a mass constraint to the loss function of the baseline PINN. This is only a soft constraint, and we call it the mass loss. It is calculated as follows:
	\begin{equation}\label{L_m}
		\mathcal{L}_m(t,\theta)=\frac{1}{N_t}\sum_{i=1}^{N_t}|m(\hat{u}_{\theta},t_m^i)-m(u,0)|^2,
	\end{equation}
	where $m(u_{\theta},t)$ denotes the mass at time $t$, which can be calculated by numerical integration methods. $\{\bm{x}_q^j,\omega_j\}_{j=1}^{N_q}$ denotes the orthogonal points and the corresponding weights in $\Omega$, respectively. The mass $m(\hat{u}_{\theta},t_m^i)$ at $t_m^i$ is calculated as:
	\begin{equation}
		m(\hat{u}_{\theta},t_m^i)=\int_{\Omega}\hat{u}_{\theta}(\bm{x},t_m^i)d\bm{x}\approx\sum_{j=0}^{N_q}\omega_j \hat{u}_{\theta}(\bm{x}_q^j,t_m^i).
	\end{equation}\par
	The above mass loss (\ref{L_m}) is then added to the PINN loss function as a regularization term. It enables the predicted solutions to maintain the initial mass and also improves the accuracy of the PINN.\par

    \subsection{Mass-preserving spatio-temporal adaptive PINN}
In this section, we present the complete procedure of the mass-preserving spatio-temporal adaptive PINN. The mass constraint (\ref{L_m}) is added to the PINN loss function as described in Section \ref{mass_conserivate}, and the adaptive sampling strategy from Section \ref{Sampling} is integrated into the training process.

As discussed in Section \ref{Sampling}, the subinterval length can be non-uniform and adjusted according to the difficulty of learning. For phase field equations, regions where the phase structure evolves rapidly are harder to learn, while slowly changing parts are easier. 
Since the evolution of the phase structure is consistent with the energy decay, we divide the full time domain $[0,T]$ into $N$ non-uniform segments
	\begin{equation}
		[T_0,T_1],[T_1,T_2],\cdots,[T_{N-1},T_N],\	
	\end{equation}
where $T_0=0$, $T_N=T$, and the partition points are determined by the rate of energy decrease using the threshold $\alpha$ defined in Section 2.2.2. Specifically, smaller time subdomains are assigned to regions where the energy drops rapidly, and larger ones to regions where it decays slowly.

An independent neural network is established in each time interval:
	\begin{equation}
		Net_1,Net_2,\cdots,Net_N.
	\end{equation}
Although they can be set to different network structures, for simplicity,  the same network structure is used in this paper. The initial condition of $Net_1$ is set to the given initial condition of the equation, and the output of $Net_n$ at $T_n$ ($n=1,\cdots,N-1$) serves as the initial condition for $Net_{n+1}$. During training, spatial adaptive resampling (Section \ref{Sampling}) is applied after a certain number of iterations to reinforce the training result. Once each subnetwork is sufficiently trained, the solutions on all subintervals are concatenated to obtain the full-time solution, denoted as $\hat{u}_{crude}$.\par

The following Algorithm~\ref{alg1} describes the proposed adaptive PINN: 
	\begin{algorithm}[H]\footnotesize                 
		\caption{Mass-preserving spatio-temporal adaptive PINN}          
		\label{alg1}      
		\begin{algorithmic}[1] 
			\State Divide the time domain equally into $N$ subintervals;
			\For{$n=1$ to  $N$ }
			\State Initialize the parameters of PINN with loss function:$\mathcal{L}_n=\lambda_{u}\mathcal{L}_u+\mathcal{L}_b+\mathcal{L}_r+\lambda_{m}\mathcal{L}_m$;
			\State Sample the initial,boundary and residual points using LHS;
			\State Train the PINN for a certain number of iterations;
			\State Adaptive resampling in space repeated several times
			\State Get the predictions $\hat{u}_\theta(x,t)$;
			\EndFor
			\State Plot the numerical energy and mass in the full time domain according to $\hat{u}_\theta(x,t)$;
			\State Calculate the mass error;
			\Repeat
			\State Adjust the time subintervals according to the rate of energy decrease, and fine-tune the weight of mass loss according to the mass error;
			\State Step $2\sim 8$;
			\Until{The total number of iterations or mass error and solution error reach the limit}
			\State	Get the final predicted solutions $\hat{u}_\theta(x,t)$;\\
			\Return $\hat{u}_\theta(x,t)$;
		\end{algorithmic}
	\end{algorithm}
	\section{Numerical results}
	In this section, several numerical results for the Cahn-Hilliard equation with different energy potentials will be presented. These results serve to evaluate the ability of the proposed modified method in solving parabolic equations, especially the equations with high-order spatial derivatives and singular energy functionals, as well as the effect of dimensionality. It will be illustrated by the solution errors and the mass error.\par
	The reference solutions $u$ for the following numerical results were obtained using classical numerical methods, employing a semi-implicit discretization in time and spectral discretization in space. The accuracy of the predicted solution $\hat{u}_\theta$ is measured by the absolute error and relative $L^2$ error, computed on the data point set $\{(\bm{x}_i,t)\}_{i=1}^N\in\Omega\times[0,T]$. \par
	The absolute error at time $t$ is: 
	\begin{equation}
		Error_{abs}(\bm{x}_i,t)=|\hat{u}_\theta(\bm{x}_i,t)-u(\bm{x}_i,t)|.
	\end{equation}
	The relative $L^2$ error at time $t$ is: 
	\begin{equation}
		Error_{L^2}=\frac{\sqrt{\sum_{i=1}^{N}|\hat{u}_\theta(\bm{x}_i,t)-u(\bm{x}_i,t)|^2}}{\sqrt{\sum_{i=1}^{N}|u(\bm{x}_i,t)|^2}}.
	\end{equation}
	The conservation of mass is measured by the mass error, which is defined as follows:
	\begin{equation}
		Error_{mass}=|m(\hat{u}_\theta(\bm{x},t))-m(u_0)|.
	\end{equation}\par
    
	In the rest of this section, we first test the fourth-order Cahn-Hilliard equation with Ginzburg-Landau potential, demonstrating the operation details and effectiveness of this method. Then, we extend its application to solve the Cahn-Hilliard equation with Flory-Huggins potential, which exhibits singularity. Finally, in order to further validate the effect of the mass-preserving spatio-temporal adaptive PINN, we apply the method to more complex phase field models by simulating the system of Cahn-Hilliard equations and a three-dimensional Cahn-Hilliard equation.

All numerical experiments in this section are implemented in Python 3.8 with TensorFlow 1.15.5, on a server equipped with an Intel(R) Xeon(R) Platinum 8255C CPU, an RTX 3080 GPU, and 40 GB RAM. For the two-dimensional benchmark problem (the Cahn-Hilliard equation with the Flory-Huggins potential), the training time of our method is approximately 10 hours. Compared with the standard PINN, our method requires about 25\% more training time, while improving the relative $L^2$ error by more than one order of magnitude and maintaining stable performance over time. For the three-dimensional case, the training time scales proportionally with the increased problem size.
	\subsection{Cahn-Hilliard equation with Ginzburg-Landau potential in two dimensions(2D)}\label{GL_2D}
	We start by testing a relatively simple case,  which is the 2D Cahn-Hilliard equation with the polynomial Ginzburg-Landau potential. In this example, we consider the equation given by
	\begin{equation}\label{CH_example}
		\left\{
		\begin{split}
			u_t&=\Delta\varphi,\quad (x,y)\in[-1,1]^2,\quad t\in[0,1],\\
			\varphi&=-\varepsilon^2\Delta u+(u^3-u),\\
		\end{split}
		\right.
	\end{equation}
	with  periodic boundary conditions, the parameter $\varepsilon=0.05$, and the initial profile for $u$ as
	\begin{equation}\label{initial_condition_CH2D}
		\begin{split}
			u(x,y,t=0)=\max\bigg(\tanh\frac{r-R_1}{2\varepsilon},\tanh\frac{r-R_2}{2\varepsilon}\bigg),
		\end{split}
	\end{equation}
	where $r=0.4,R_1=\sqrt{(x-0.7r)^2+y^2},R_2=\sqrt{(x+0.7r)^2+y^2}$. This initial condition consists of two circles which will simulate the process of their merging. First, we consider the case without a mass loss term, and choose the uniform size $\Delta t=0.2$ to uniformly divide the full time domain into five segments. Feedforward fully connected neural networks are established in each subinterval separately, and each network has 10 hidden layers with 50 neurons per layer, and two outputs. The hyperbolic tangent function is chosen for the activation function, and the loss function for each network is defined as follows:
	\begin{equation}
		\left\{
		\begin{split}
			\mathcal{L}=&\lambda_u \mathcal{L}_u+\mathcal{L}_b+\mathcal{L}_r,\\
			\mathcal{L}_u=&\frac{1}{N_u}\sum_{i=1}^{N_u}|\hat{u}_\theta(\bm{x}_u^i,0)-u_i|^2,\\
			\mathcal{L}_b=&\frac{1}{N_b}\sum_{i=1}^{N_b}(|\hat{u}_\theta(\bm{x}_u^i,t_b^i)-\hat{u}_\theta(\bm{x}_l^i,t_b^i)|^2+\bigg|\frac{\partial\hat{u}_\theta }{\partial \bm{x}}(\bm{x}_u^i,t_b^i)-\frac{\partial\hat{u}_\theta}{\partial \bm{x}}(\bm{x}_l^i,t_b^i)\bigg|^2),\\
			\mathcal{L}_r=&\frac{1}{N_r}\sum_{i=1}^{N_r}(|\hat{\varphi}_\theta(\bm{x}_r^i,t_r^i)+\varepsilon^2\Delta \hat{u}_\theta(\bm{x}_r^i,t_r^i)-(\hat{u}_\theta(\bm{x}_r^i,t_r^i)^3-\hat{u}_\theta(\bm{x}_r^i,t_r^i))|^2\\ 
			&+\bigg|\frac{\partial\hat{u}_\theta}{\partial t}(\bm{x}_r^i,t_r^i)-\Delta\hat{\varphi}_\theta(\bm{x}_r^i,t_r^i)\bigg|^2).
		\end{split}
		\right.
	\end{equation}\par
	In this experiment, 10,000 initial points, 1,600 boundary points, and 20,000 residual points are used for training. The weight $\lambda_u$ of the initial condition is set to 100, and all other weights are $1$. To train the neural network, we first use the Adam optimizer with a learning rate of 0.001, and then the L-BFGS-B optimizer to fine-tune the neural network. \par
	Adaptive resampling in space is performed every 5,000 steps of training, and each time 1000 residual points with large residual value are selected to be added to the initial set of residual points. After experimental testing, we chose to perform three times spatial adaptive resampling because the accuracy and computational cost are optimal in this case.\par
\begin{figure}[H]
		\centering         
		\subfloat[Predicted solutions at $t=0,0.25,0.5,1$]   
		{
			\includegraphics[width=1\textwidth]{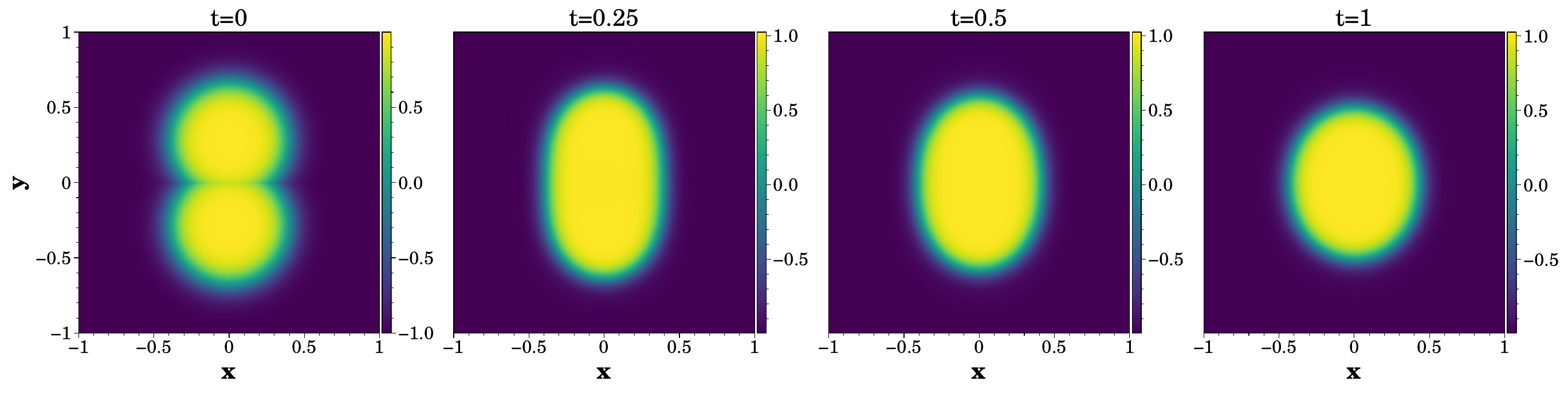}
			\label{result_1}
		}\\
		\subfloat[Absolute errors at $t=0,0.25,0.5,1$]
		{
			\includegraphics[width=1\textwidth]{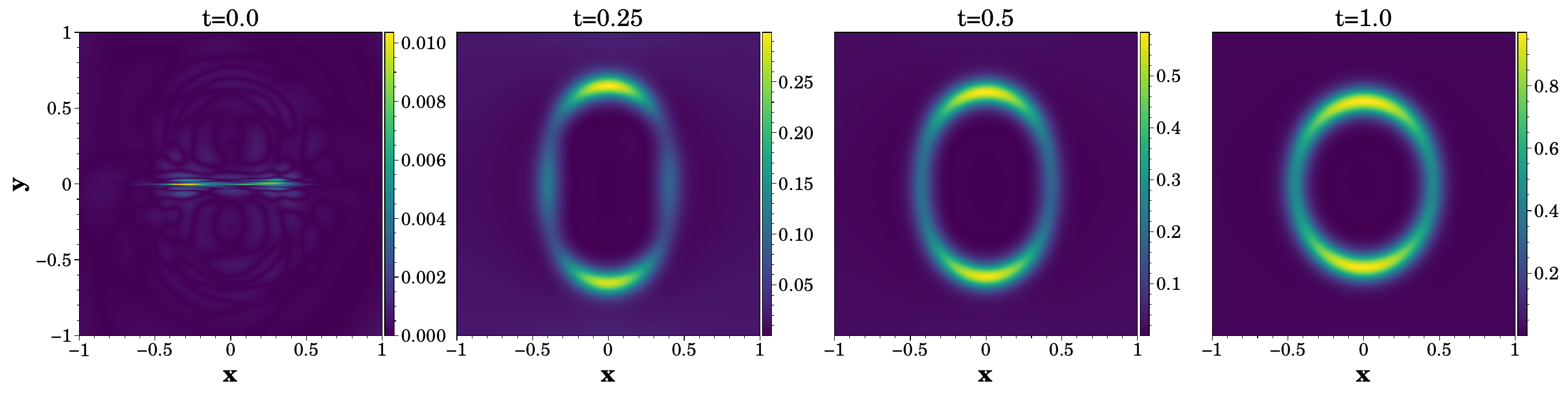}
			\label{error_1}
		}
		\caption{Predicted solutions and absolute errors of 2D Cahn-Hilliard equation with Ginzburg-Landau potential given by spatio-temporal adaptive PINN with $\Delta t=0.2$.}\label{1_RE}
	\end{figure}

	The predicted solutions at different times ($t=0,0.25,0.5,1$) and the corresponding absolute errors with uniform $\Delta t=0.2$ are shown in Fig.~\ref{1_RE}\subref{result_1} and  Fig.~\ref{1_RE}\subref{error_1}. 
    Although the merging of two circles is captured, the maximum absolute errors of the predicted solutions are large, with values of $0.298, 0.583$ and $0.971$, and increase with time. Thus, further improvement is needed to reduce the absolute errors.\par

\begin{figure}[H]
		\centering   
		\setcounter{subfigure}{0}      
		\subfloat[Energy]   
		{
			\includegraphics[width=0.45\textwidth]{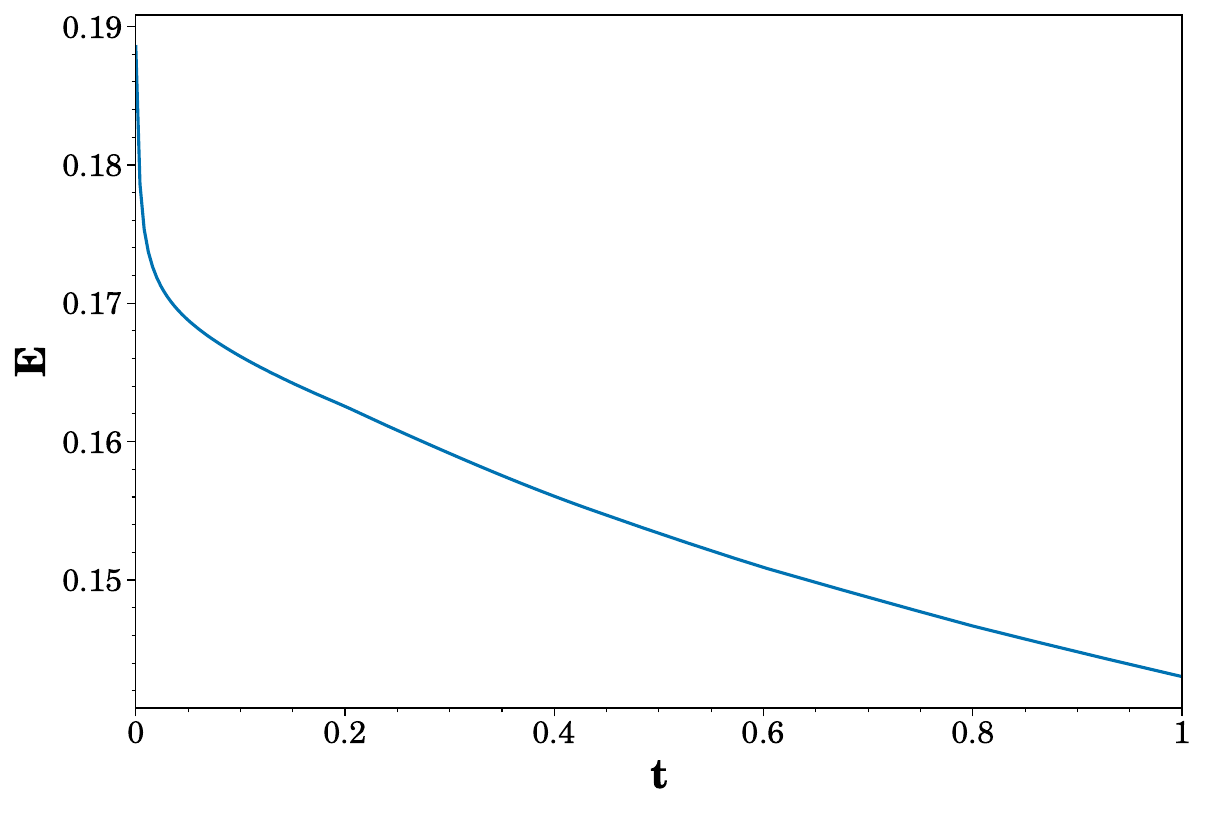}
			\label{energy_1}
		}
		\subfloat[Mass]
		{
			\includegraphics[width=0.46\textwidth]{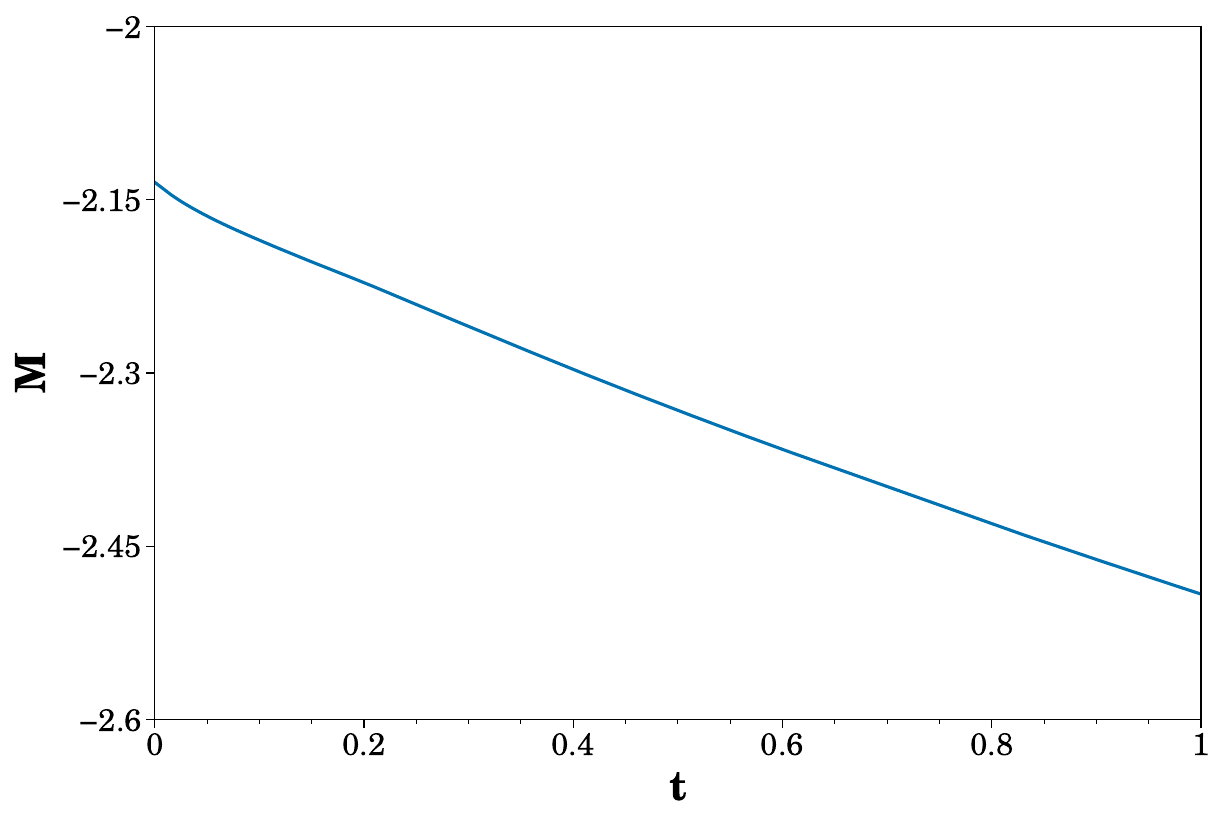}
			\label{mass_1}
		}
		\caption{The numerical energy(E) and mass(M) of the predicted solutions $\hat{u}_\theta(\bm{x},t)$.}
		\label{1_EM}
	\end{figure}

	The energy and mass variations are shown in Fig.~\ref{1_EM}\subref{energy_1} and  Fig.~\ref{1_EM}\subref{mass_1}. The mass is not conserved although the energy decreases. The shrinking yellow region in Fig.~\ref{1_RE}\subref{result_1} also indicates mass loss. The original system has the property of mass conservation, but the PINN predictions do not comply with this principle. Therefore, it is necessary to make adjustments to ensure that the predictions of the PINN adhere to physical laws.
	\par
 
    \vspace{-0.3cm}
	\par	
	\vspace{-0.5cm}
	Therefore, we add the following mass loss term $\mathcal{L}_m$ to the total loss function for each PINN:
	\begin{equation}\label{loss_m}\begin{split}
			&\mathcal{L}_m=\frac{1}{N_t}\sum_{i=1}^{N_t}|m(\hat{u}_\theta,t_m^i)-m(u,0)|^2,\\
		\end{split}
	\end{equation}
	where
	\begin{equation}
		m(u,t)=\int_{\Omega}ud\bm{x}\approx\frac{1}{N_q}\sum_{j=1}^{N_q}u(\bm{x}_q^j,t).
	\end{equation}\par

\begin{figure}[H]
		\centering        
		\subfloat[Predicted solutions at $t=0,0.25,0.5,1$]   
		{
			\includegraphics[width=1\textwidth]{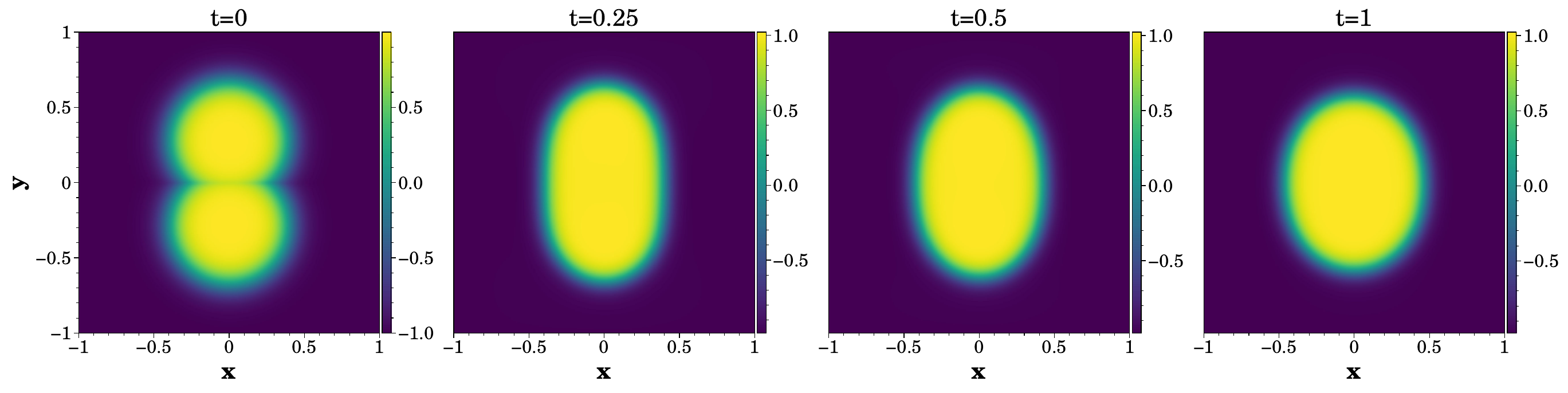}
			\label{result_2}
		}\\
		\subfloat[Absolute errors at $t=0,0.25,0.5,1$]
		{
			\includegraphics[width=1\textwidth]{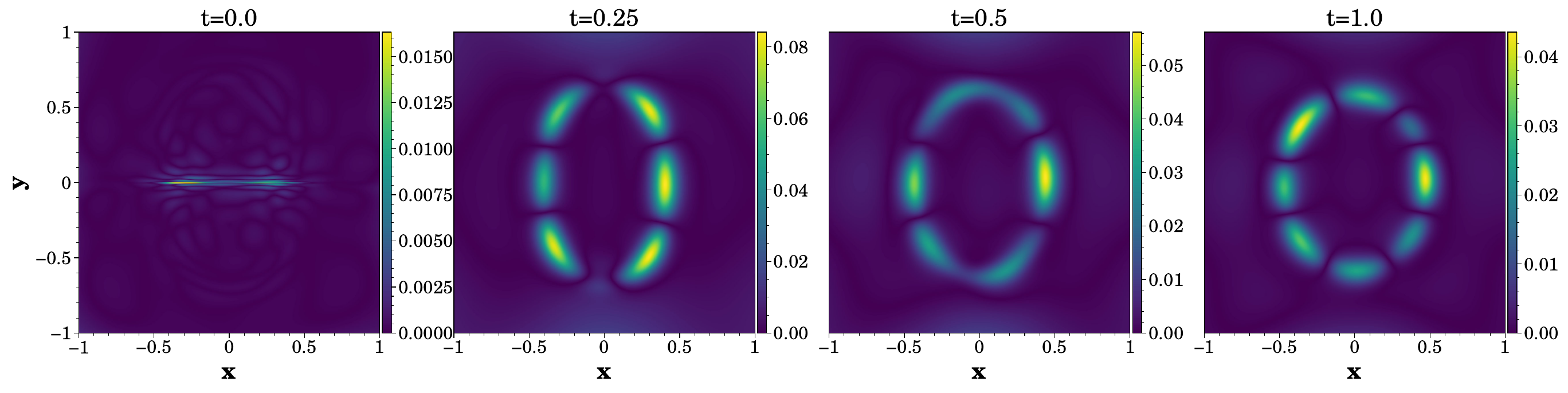}
			\label{error_2}
		}
		\caption{Predicted solutions and absolute errors of 2D Cahn-Hilliard equation with Ginzburg-Landau potential, with mass loss and adaptive sampling.}
		\label{2_RE}
	\end{figure}

	\vspace{-0.3cm}	
	We choose $N_t=21$ equidistant points in each time subinterval as $t_m^i$ in Eq.(\ref{loss_m}), and the mass at time $t_i$ is computed at $50\times 50$ quadrature points $\{(\bm{x}_q^j,y_q^j)\}_{j=1}^{N_q}$.
	The predicted solutions and absolute errors are shown in Fig.~\ref{2_RE}\subref{result_2} and  Fig.~\ref{2_RE}\subref{error_2}. From Fig.~\ref{2_RE}\subref{error_2}, it is observed that the maximum absolute error at each time decreases to $0.084, 0.056$ and $0.0444$, which are smaller than the absolute errors in Fig.~\ref{1_RE}\subref{error_1}.
	
	\vspace{-0.3cm}	
	\par 
	\vspace{-0.5cm}

\begin{figure}[H]
		\centering
		\subfloat[Energy]   
		{
			\includegraphics[width=0.45\textwidth]{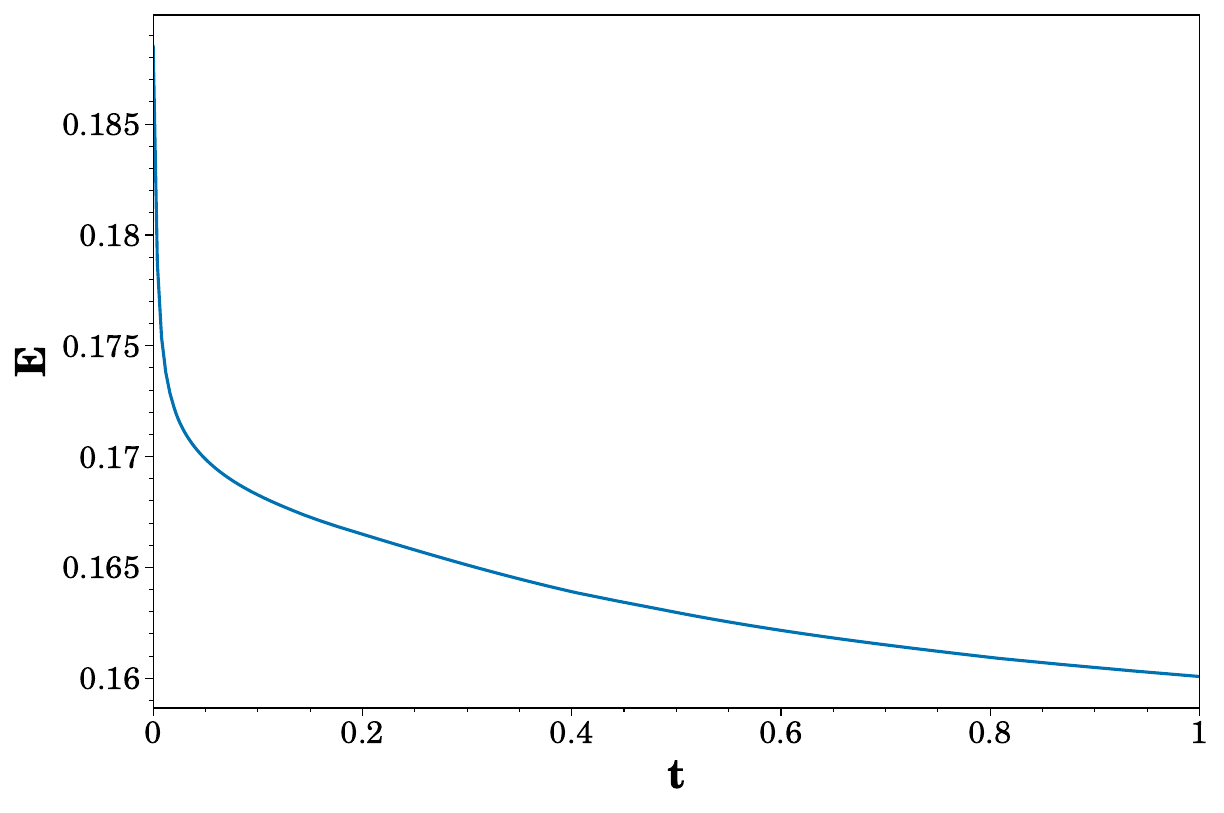}
			\label{energy_2}
		}
		\subfloat[Mass]
		{
			\includegraphics[width=0.45\textwidth]{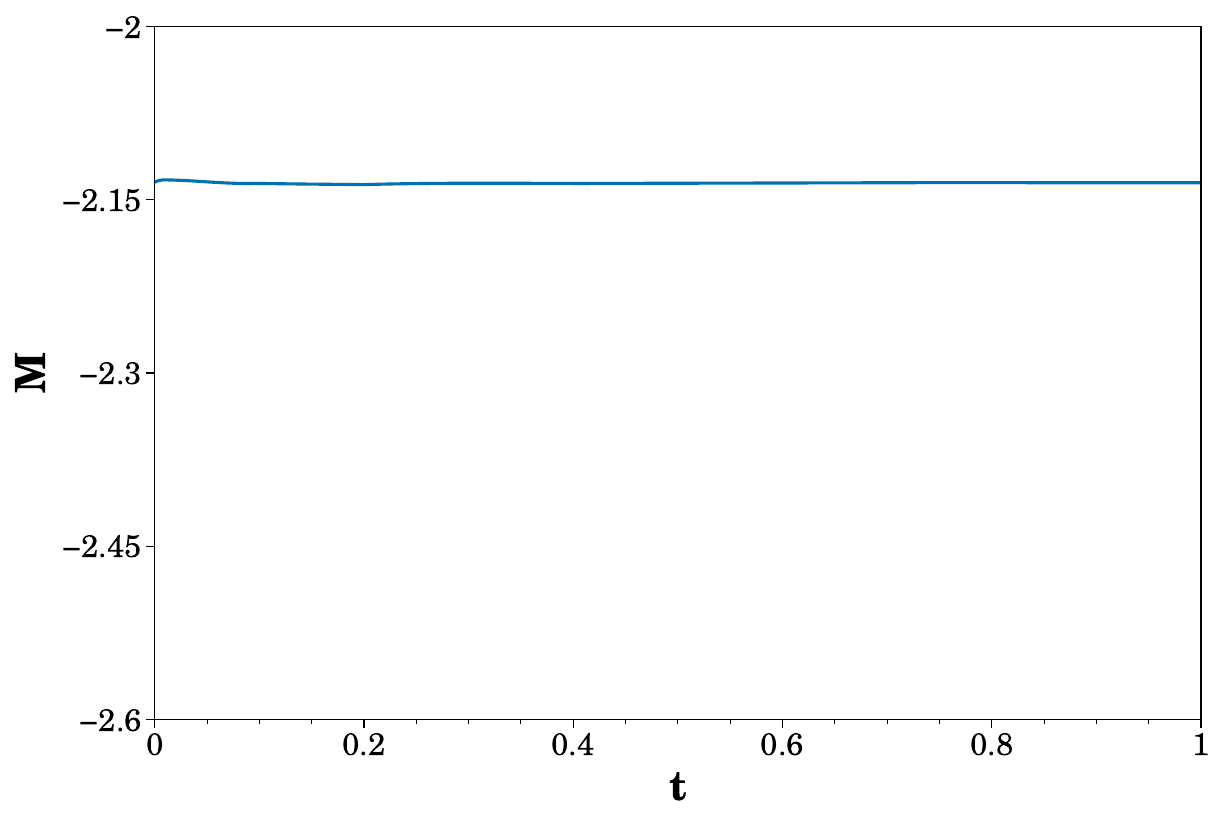}
			\label{mass_2}
		}
		\caption{The numerical energy(E) and mass(M) of the predicted solutions $\hat{u}_\theta(\bm{x},t)$.}
		\label{2_EM} 
	\end{figure}

After adding the mass loss function, the energy and mass evolutions are plotted in Fig.~\ref{2_EM}\subref{energy_2} and Fig.~\ref{2_EM}\subref{mass_2}. The energy of the predicted solutions decreases with time. The mass is relatively stable in the later stages, while the initial mass error is large, with a maximum of $2.27e-03$ over the entire time domain.\par

\begin{figure}[H]
		\centering         
		\subfloat[Predicted solutions at $t=0,0.25,0.5,1$]   
		{
			\includegraphics[width=1\textwidth]{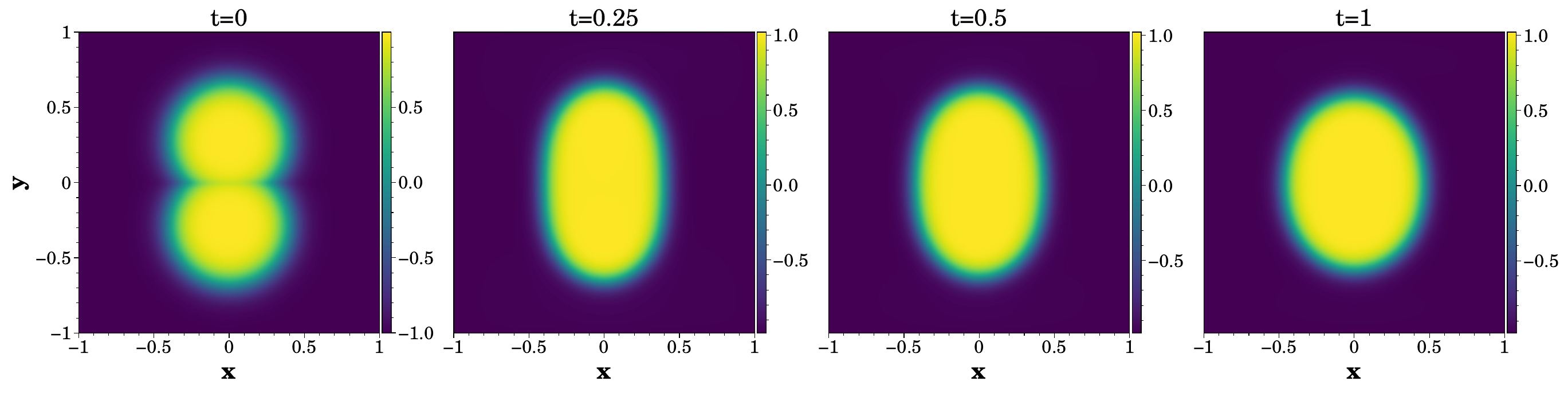}
			\label{result_4}
		}\\
		\label{4}
		\subfloat[Absolute errors at $t=0,0.25,0.5,1$]
		{
			\includegraphics[width=1\textwidth]{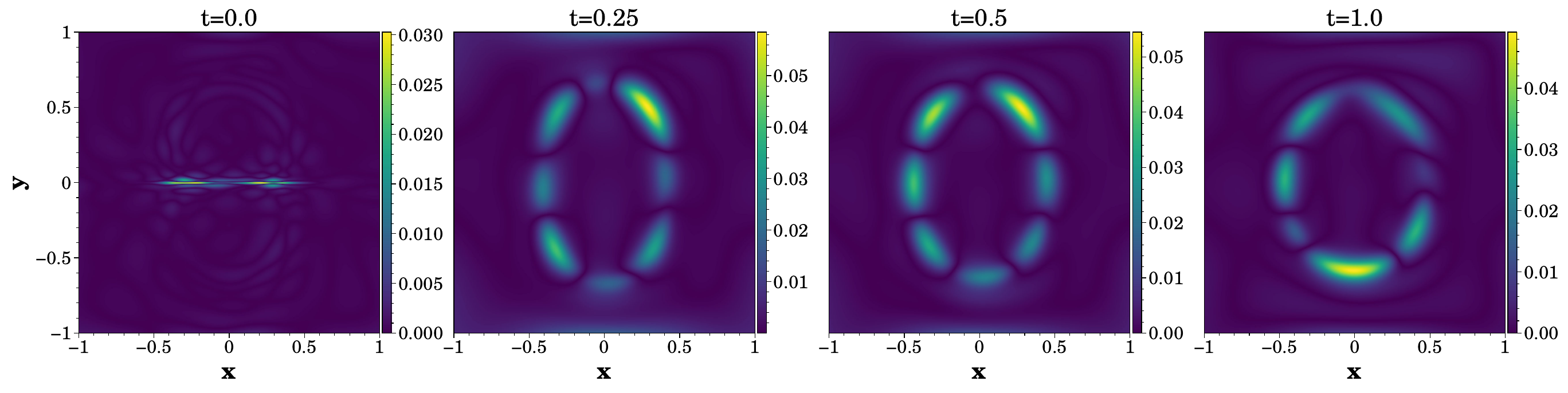}
			\label{error_4}
		}
		\caption{Predicted solutions and absolute errors of the 2D Cahn-Hilliard equation with Ginzburg-Landau potential, with time subintervals of lengths $0.02, 0.18, 0.4, 0.4$.} 
		\label{4_RE}
	\end{figure}

    To reduce the initial mass error, time subintervals are adjusted according to the variation of energy with time. The energy decreases rapidly in $t\in[0,0.02]$ and then slowly, becoming nearly constant after $t=0.2$. Therefore, we reset the lengths of time subintervals to $0.02, 0.18, 0.4, 0.4$, with $N_t=26$. The mass loss weight in the first time subdomain is set to $10$, and the weights in the other segments are kept as $1$. The predicted solutions and absolute errors are shown in Fig.~\ref{4_RE}\subref{result_4} and Fig.~\ref{4_RE}\subref{error_4}. The maximum absolute errors at each time in Fig.~\ref{4_RE}\subref{error_4} are $0.058, 0.054$ and $0.049$. Although not improved significantly compared to the
     uniform case, the maximum mass error is reduced to $1.01e-03$. \par
	\vspace{-0.5cm}

	\par
\begin{figure}[H]
		\centering
		\setcounter{subfigure}{0}  
		\subfloat[Energy]   
		{
			\includegraphics[width=0.45\textwidth]{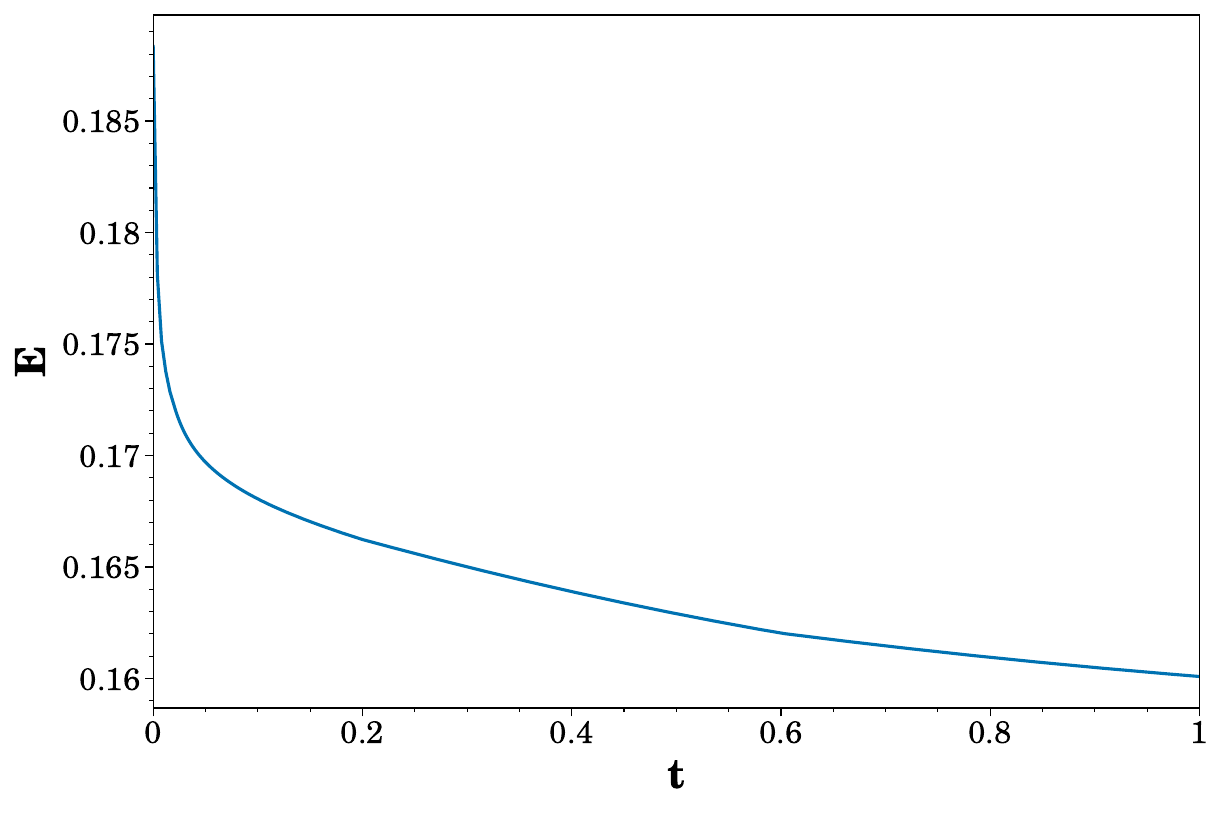}
			\label{energy_4}
		}
		\subfloat[Mass]
		{
			\includegraphics[width=0.45\textwidth]{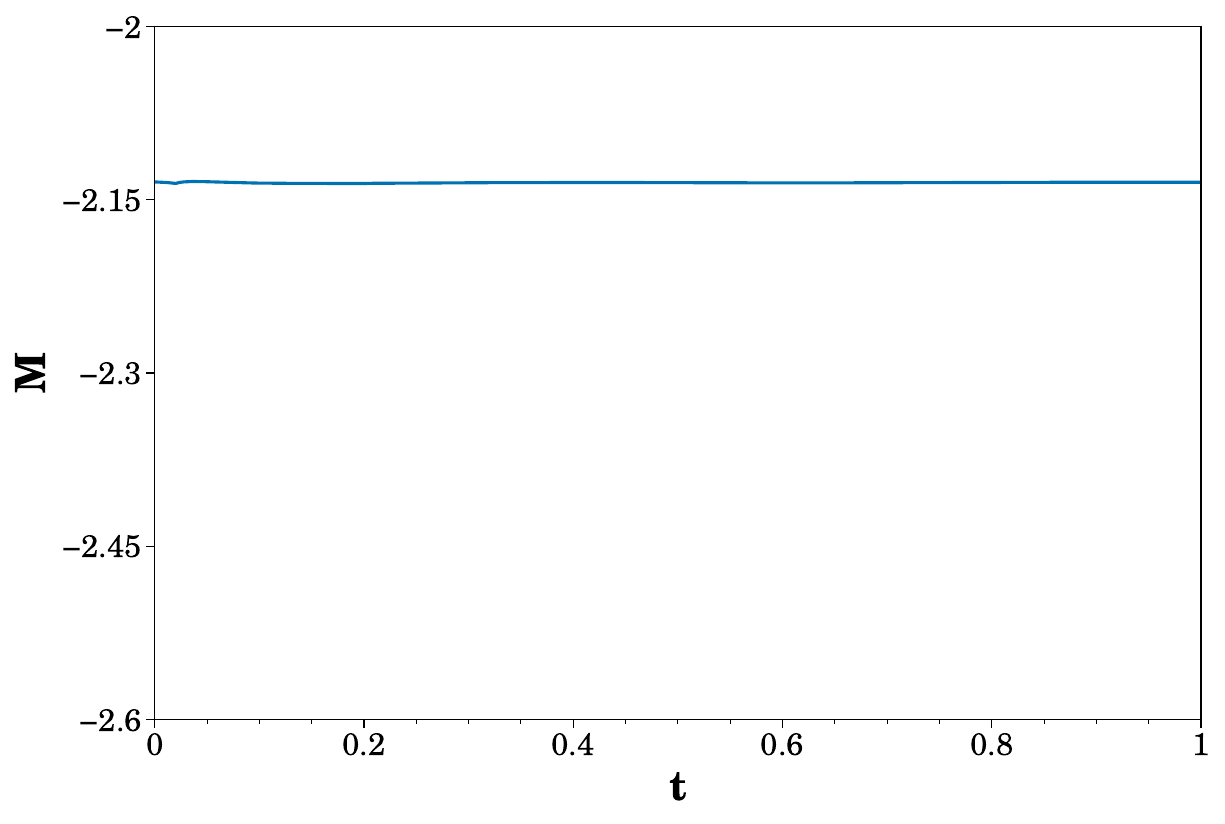}
			\label{mass_4}
		}
		\caption{The numerical energy(E) and mass(M) of the predicted solutions $\hat{u}_\theta(\bm{x},t)$.}
		\label{4_EM}
	\end{figure}

 
	\vspace{-0.5cm}
	The numerical energy and mass are shown in Fig.~\ref{4_EM}\subref{energy_4} and Fig.~\ref{4_EM}\subref{mass_4}. Fig.~\ref{4_EM}\subref{mass_4} shows that the total mass keeps nearly constant throughout the simulation..\par
	To better understand the contribution of each component, we conduct an ablation study by comparing three variants. A comparison of the relative $L^2$ errors and maximum mass errors is shown in Table~\ref{table1}. In Table~\ref{table1}, $err_1$ represents the errors obtained by solving the equation (\ref{CH_example}) only by using adaptive sampling in time and space with uniform size. $err_2$ are errors obtained by adding a mass constraint based on the $err_1$ operation. $err_3$ corresponds to the full model, which further adjusts the time subintervals to be non-uniform (according to the rate of energy decrease) and fine-tunes the mass loss weights. We change the mass loss weight to $10$ for the first time subinterval and keep the weights for the others as 1. As shown in Figs.~4(b), 6(b), and 8(b), the absolute errors at $t=0.25$ and $t=0.5$ decrease progressively as each component is added. The full model achieves the lowest relative $L^2$ error at $t=0.25$ and $t=1$ and the smallest maximum mass error, confirming the effectiveness of each proposed component. \par
	\begin{table}[H]
		\caption{Cahn-Hilliard equation with Ginzburg-Landau potential: ablation study -- relative $L^2$ errors and maximum mass error of three cases at different times}
		\label{table1}
		\centering
		\def\temptablewidth{0.8\textwidth}
			\begin{tabular*}{\temptablewidth}{@{\extracolsep{\fill}}llcccc}
				\hline
				&Time&$err_1$&$err_2$&$err_3$&\\ 
				\hline
				&0   &5.23e-04&7.26e-04&1.39e-03&\\
				&0.25&5.10e-02&1.59e-02&8.36e-03&\\
				&0.5 &1.12e-01&8.34e-03&8.40e-03&\\
				&1   &2.15e-01&7.13e-03&7.10e-03&\\
				\hline
				&mass error&3.56e-01&2.27e-03&1.01e-03&\\
				\hline
			\end{tabular*}\\
		\end{table}

We note that such a process can be applied to a large class of complex parabolic equations, especially to phase field equations. Although the error reduction by this procedure is not significant for the case of Ginzburg-Landau potential, the improvement becomes significant for the Flory-Huggins potential. It significantly reduces the absolute errors of the predicted solutions and helps preserve mass, as will be shown in the following.
        
    \begin{table}[H]
            \caption{Cahn-Hilliard equation with Ginzburg-Landau potential: sensitivity analysis of the energy threshold --- relative $L^2$ errors and maximum mass errors for two energy-threshold choices and the final adaptive setting.} 
			\label{tab:sensitivity}
			\centering
			\def\temptablewidth{0.8\textwidth}
				\begin{tabular*}{\temptablewidth}{@{\extracolsep{\fill}}llcccc}
					\hline
					&Time&$\alpha=10\%$&$\alpha=5\%$&Final&\\ 
					\hline
					&0   &1.24e-03&1.57e-03&1.39e-03 &\\
					&0.25&3.85e-02&8.54e-03&8.36e-03 &\\
					&0.5 &8.11e-02&4.57e-03&8.40e-03 &\\
					&1   &1.52e-01&1.00e-02&7.10e-03 &\\
					\hline
					&mass error&9.02e-04&6.84e-03&1.01e-03\\
					\hline
				\end{tabular*}
		\end{table}

    We further perform a sensitivity analysis on the energy threshold parameter $\alpha$ (defined in Section 2.2.2). The results, summarized in Table \ref{tab:sensitivity}, show that the final threshold achieves the best balance between long-time accuracy and mass conservation, confirming the reasonableness of our configuration. We acknowledge that the current choice is heuristic and discuss the possibility of a more theory-driven rule as future work.

		\subsection{Cahn-Hilliard equation with Flory-Huggins potential in two dimensions(2D)}\label{FH_2D}
		
		In this section, we test the 2D Cahn-Hilliard equation with Flory-Huggins potential. It has a highly nonlinear and singular energy potential, so we need to reduce the error due to singularity. The computational domain is $\Omega:=[-1,1]^2$, $t\in[0,1]$. We focus on
		the following specific form,
		\begin{equation}\label{FH_example}
			\left\{
			\begin{split}
				u_t&=\Delta\varphi,\\
				\varphi&=-\varepsilon^2\Delta u+(\ln(1+u)-\ln(1-u)-\theta u),\\
			\end{split}
			\right.
		\end{equation}
		with  periodic boundary conditions. We set $\varepsilon=0.05$, $\theta=3$, along with the initial condition
		\begin{equation}\label{initial_condition_FH2D}
			\begin{split}
				u(x,y,t=0)=\max\bigg(\tanh\frac{r-R_1}{2\varepsilon},\tanh\frac{r-R_2}{2\varepsilon}\bigg)\times 0.9,
			\end{split}
		\end{equation}
		where $r=0.4,R_1=\sqrt{(x-0.7r)^2+y^2},R_2=\sqrt{(x+0.7r)^2+y^2}$.\par

		The Ginzburg-Landau potential energy is in polynomial form and does not require special handling in computations. But the Flory-Huggins potential contains logarithmic terms, and the result will blow up as $u$ approaches $\pm 1$ due to its singularity. Therefore, we usually truncate $u$ at each time subinterval to ensure that its value does not approach $\pm1 $ when using traditional numerical methods. When using PINN to solve problems with singularity, gradient explosion occurs during network training. This is mainly caused by the logarithmic term in the residual loss of the equation. 
		Therefore, we add a mapping to the output of the network such that the value of the output u is restricted to $[-0.9,0.9]$ to avoid the blow up produced by the logarithmic term. The mapping is represented as follows:
		\begin{equation}
			u_{map}=1.8\times sigmoid(u)-0.9.
		\end{equation}\par
\begin{figure}[H]
			\centering         
			\subfloat[Predicted solutions at $t=0,0.25,0.5,1$]   
			{
				\includegraphics[width=1\textwidth]{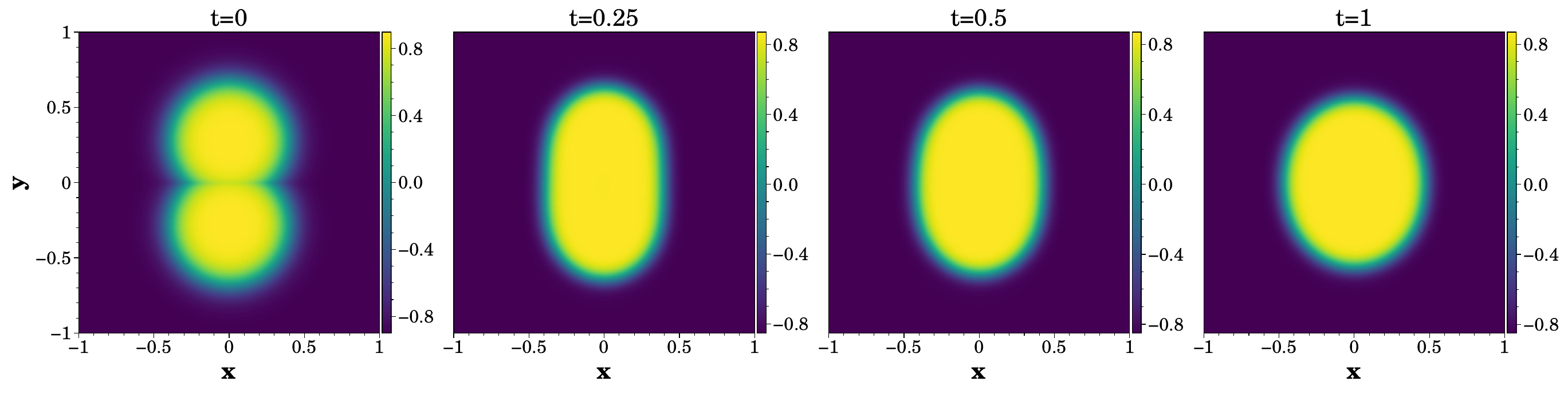}
				\label{result_5}
			}\\
			\subfloat[Absolute errors at $t=0,0.25,0.5,1$]
			{
				\includegraphics[width=1\textwidth]{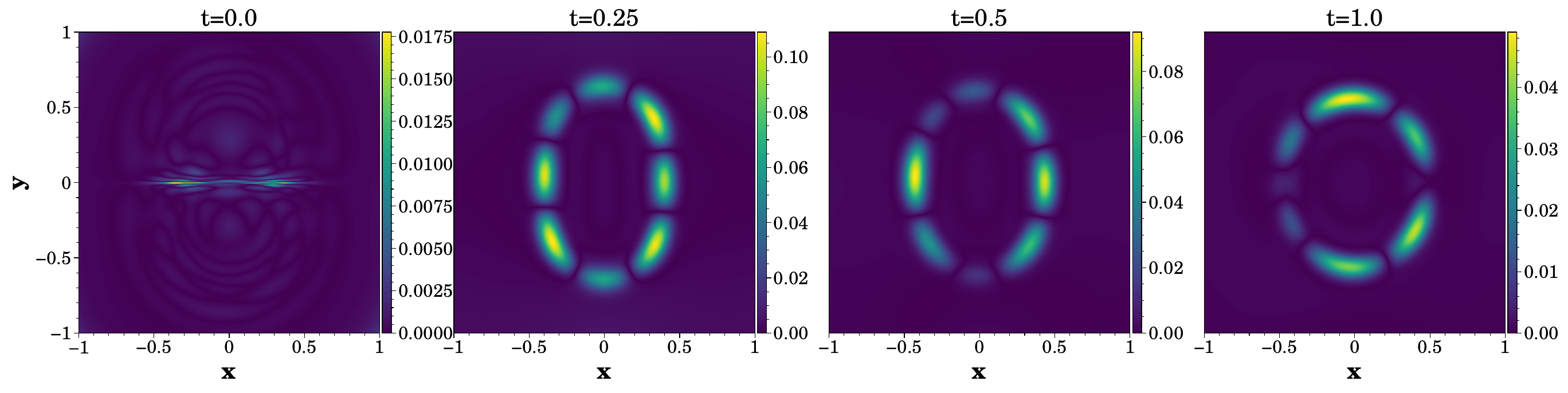}
				\label{error_5}
			}
			\caption{Predicted solutions and absolute errors of 2D Cahn-Hilliard equation with Flory-Huggins potential obtained by adaptive PINN with mass loss and spatio-temporal adaptive sampling, using uniform $\Delta t=0.2$.}
			\label{5_RE}
		\end{figure}

First, we use the same network hyperparameters as in Section~\ref{GL_2D} and divide the time domain equally into five segments, with five sub-networks established. Each sub-network performs $3$ resampling iterations with $N_t=21$. The predicted solutions and absolute errors are shown in Fig.~\ref{5_RE}\subref{result_5} and Fig.~\ref{5_RE}\subref{error_5}. The maximum absolute errors at $t=0.25,0.5,1$ are $0.109, 0.092$, and $0.049$, respectively, much larger than those in Section~\ref{GL_2D} (also with uniform segments).
		\par
\begin{figure}[H]
			\centering
			\subfloat[Energy]   
			{
				\includegraphics[width=0.45\textwidth]{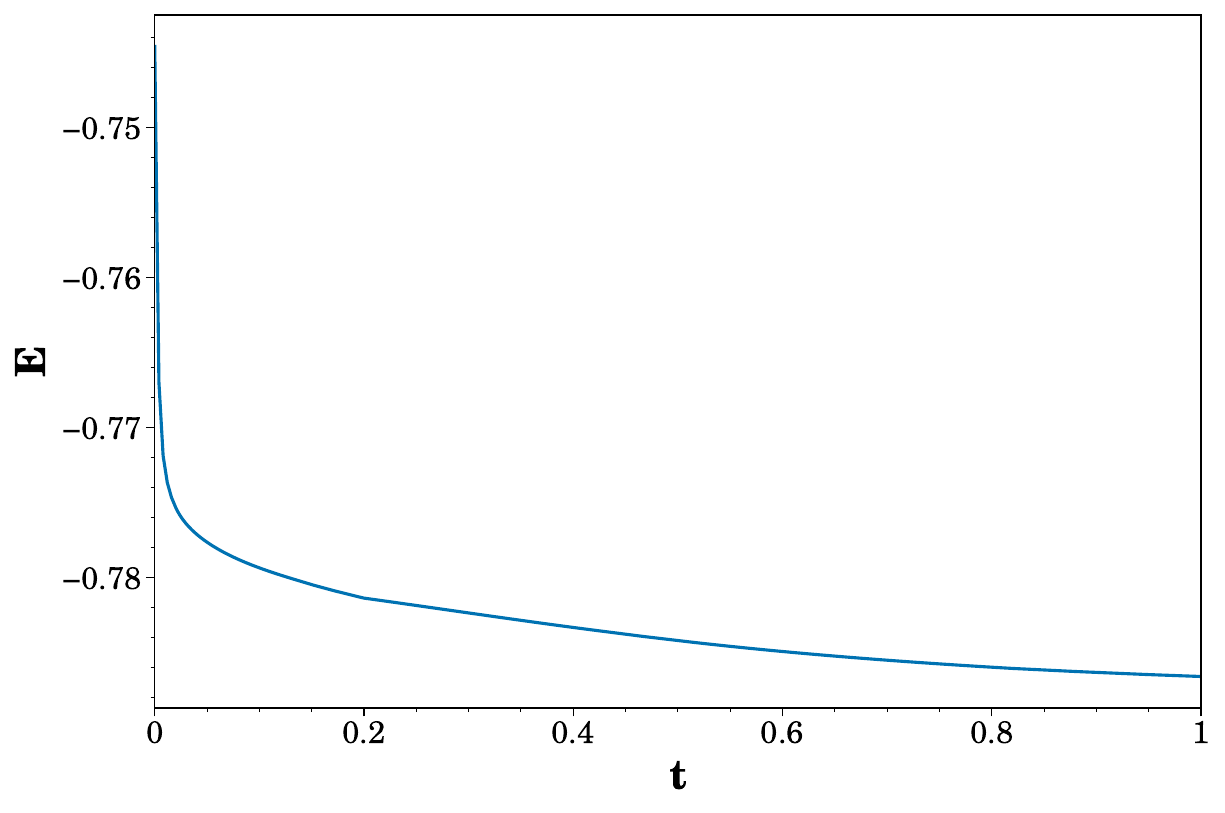}
				\label{energy_5}
			}
			\subfloat[Mass]
			{
				\includegraphics[width=0.45\textwidth]{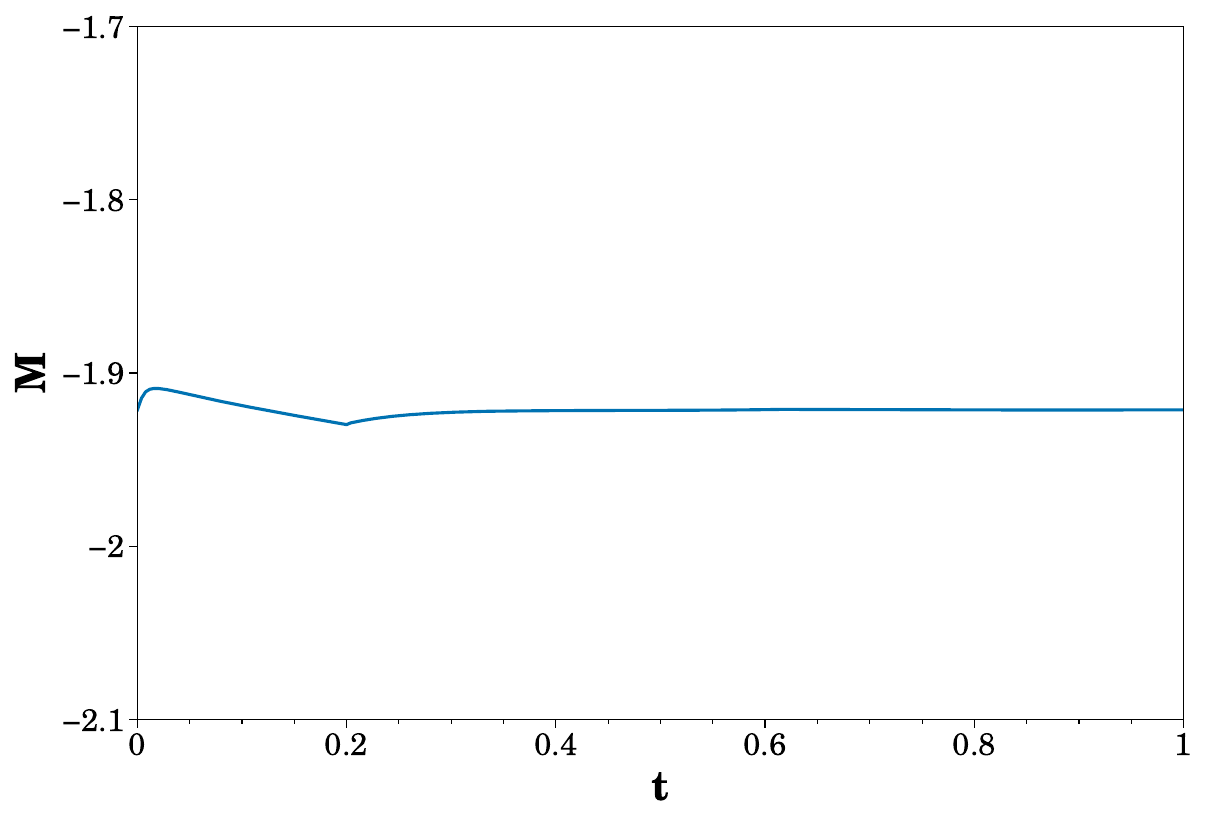}
				\label{mass_5}
			}
			\caption{The numerical energy(E) and mass(M) of the predicted solutions $\hat{u}_\theta(\bm{x},t)$.}
			\label{5_EM}
		\end{figure}

 
Based on the predicted solutions above, we plot the numerical energy and mass in Fig.~\ref{5_EM}\subref{energy_5} and Fig.~\ref{5_EM}\subref{mass_5}. The mass fluctuates significantly in $t\in[0,0.4]$ with a jump at the junction of the two intervals, and increasing the mass loss parameter $N_t$ does not help. 
As shown in Fig.~\ref{5_EM}\subref{energy_5}, the energy drops sharply in $t\in[0,0.05]$, where the mass also fluctuates severely, indicating that this period is difficult for the network to learn; the energy decreases more slowly in $t\in[0.05,1]$, which is relatively easy. Thus, we adjust the subinterval lengths according to the energy decrease rate.
\vspace{-0.3cm}

		\vspace{-0.3cm}
		\begin{figure}[H]
			\centering         
			\subfloat[Predicted solutions at $t=0,0.25,0.5,1$]   
			{
				\includegraphics[width=1\textwidth]{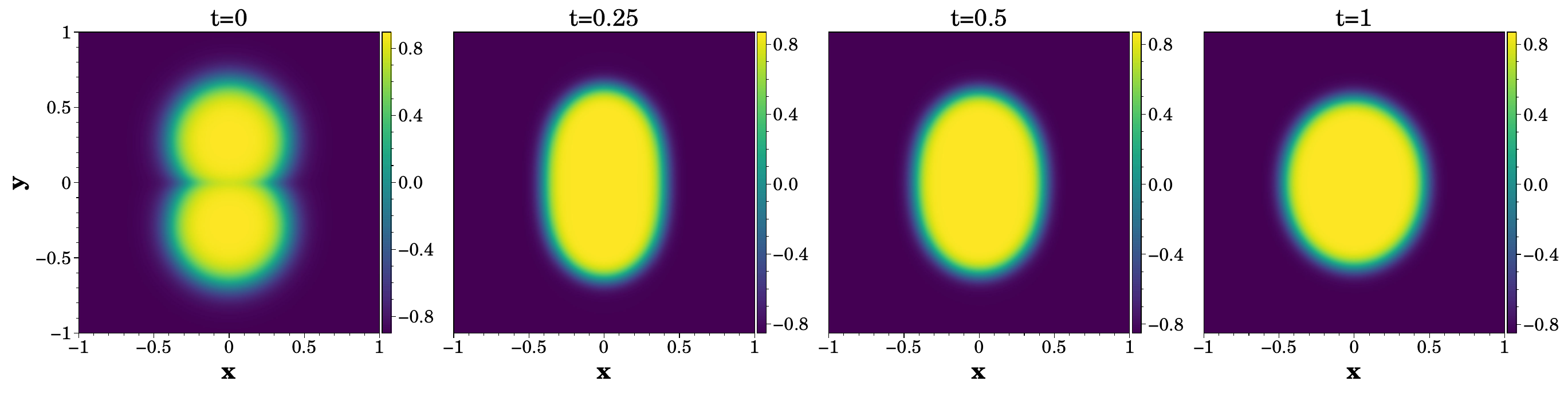}
				\label{result_7}
			}\\
			\subfloat[Absolute error at $t=0,0.25,0.5,1$]
			{
				\includegraphics[width=1\textwidth]{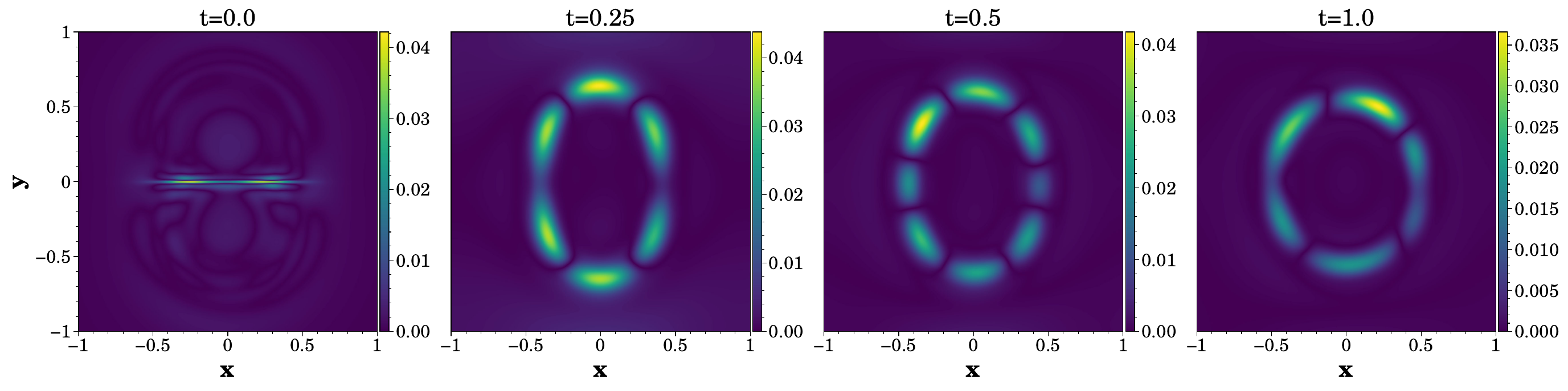}
				\label{error_7}
			}
			\caption{Predicted solutions and absolute errors of 2D Cahn-Hilliard equation with Flory-Huggins potential obtained by adaptive PINN. The mass loss weight of PINN in each interval is changed to 50, 10, 10, and 10.}
			\label{7_RE}
		\end{figure}\par 

We choose the time subintervals of lengths $0.05$, $0.15$, $0.4$ and $0.4$, with $N_t=26$ for the mass loss. The weights of mass loss for each corresponding PINNs are set to $50$, $10$, $10$, $10$, respectively. The predicted solutions and the absolute errors are shown in Fig.~\ref{7_RE}\subref{result_7} and Fig.~\ref{7_RE}\subref{error_7}. 
The maximum absolute errors of the predicted solutions at different times $(t=0.25,0.5,1)$ are $0.044, 0.042$ and $0.037$, which are greatly reduced compared to the simulation results based on the uniform subintervals.\par
        
\begin{figure}[H]
			\centering
			\subfloat[Energy]   
			{
				\includegraphics[width=0.45\textwidth]{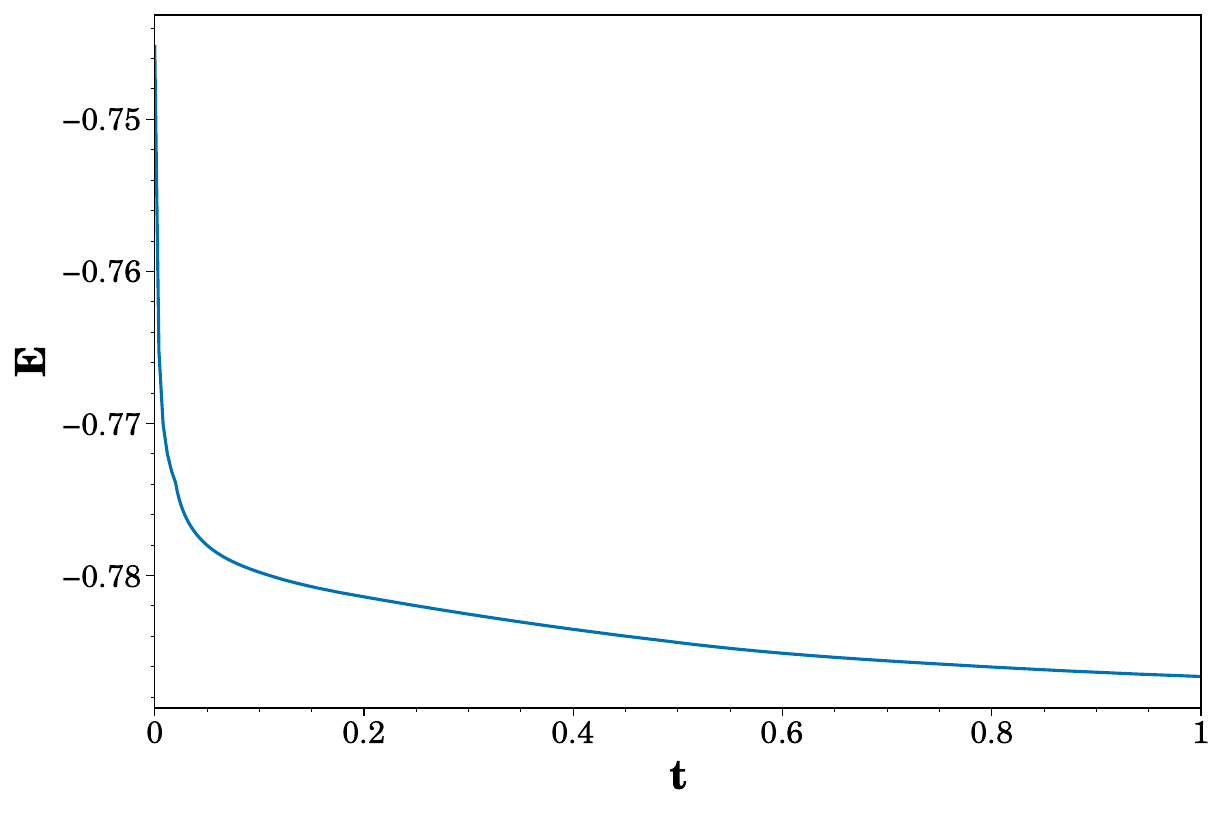}
				\label{energy_7}
			}
			\subfloat[Mass]
			{
				\includegraphics[width=0.45\textwidth]{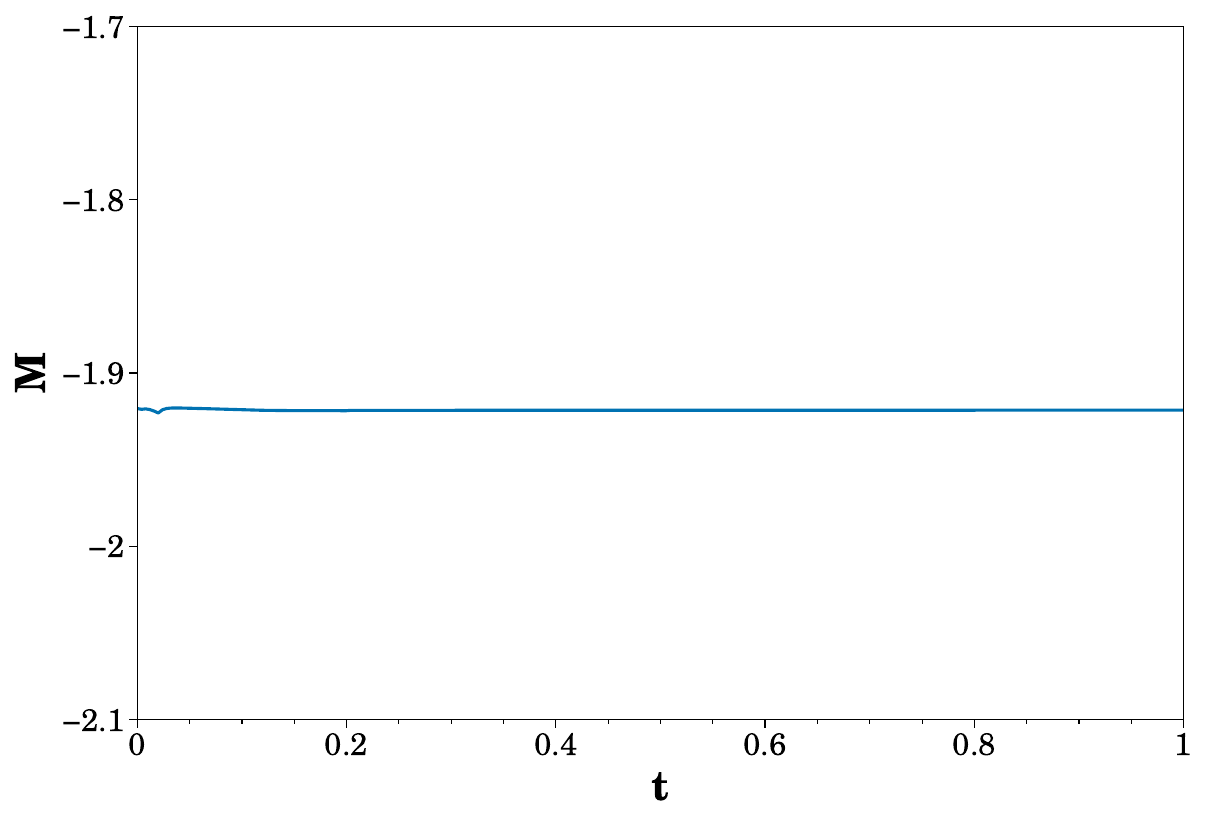}
				\label{mass_7}
			}
			\caption{The numerical energy(E) and mass(M) of the predicted solutions $\hat{u}_\theta(\bm{x},t)$.}
			\label{7_EM}
		\end{figure}

		The numerical energy and mass are plotted in Fig.~\ref{7_EM}\subref{energy_7} and Fig.~\ref{7_EM}\subref{mass_7}. The maximum mass error decreases to 1.65e-03, and the maximum absolute errors of the predictions are greatly reduced. From the above experimental results, it can be seen that the computational accuracy of the predicted solution at initial time is important for long time simulations.
		A comparison of the relative $L^2$ errors and maximum mass errors between the above two cases is shown in Table~\ref{table2}. In $err_1$, mass loss is added and the adaptive sampling method is used, with uniform segments. $err_2$ denotes the errors obtained by non-uniform segments.
		\begin{table}[H]
			\caption{Cahn-Hilliard equation with Flory-Huggins potential: relative $L^2$ errors and maximum mass errors of two cases at different times}
			\label{table2}
			\centering
			\def\temptablewidth{0.8\textwidth}
			\begin{threeparttable} 
				\begin{tabular*}{\temptablewidth}{@{\extracolsep{\fill}}llccc}
					\hline
					&Time&$err_1$&$err_2$ &\\ 
					\hline
					&0   &9.38e-04&2.69e-03 &\\
					&0.25&2.24e-02&9.30e-03 &\\
					&0.5 &1.50e-02&6.97e-03 &\\
					&1   &8.39e-03&5.75e-03 &\\
					\hline
					&mass error&1.25e-02&1.65e-03&\\
					\hline
				\end{tabular*}
			\end{threeparttable} 
		\end{table}
		Table~\ref{table2} shows that the relative $L^2$ error  at $t=0.25$, $0.5$, and $1$, and the maximum mass error are reduced when using non-uniform time partition. For the Cahn-Hilliard equation with Flory-Huggins potential, this method is significantly effective. It greatly reduces the absolute errors of the predicted solutions and helps preserve mass.
		
		\subsection{Simulations for the high-dimensional Cahn-Hilliard equation and the system of Cahn-Hilliard equations}\label{CH_3D}
		\subsubsection{Cahn-Hilliard equation with Flory-Huggins potential in three dimensions(3D)}\label{FH_3D}
		In this section, we test 3D Cahn-Hilliard equation with Flory-Huggins potential to check the ability of the algorithm to simulate high dimensional problems. The solution domain for this numerical example is $\Omega:=[-1,1]^3$, $t\in[0,1]$. The system is summarized as follows:
		\begin{equation}\label{FH_example_3D}
			\left\{
			\begin{split}
				&u_t=\Delta\varphi,\\
				&\varphi=-\varepsilon^2\Delta u+(\ln(1+u)-\ln(1-u)-\theta u).\\
			\end{split}
			\right.
		\end{equation}
		We choose the periodic boundary conditions, and set the parameters as $\varepsilon=0.05$, $\theta=3$. The initial condition consists of two spheres,
		\begin{equation}\label{initial_condition_FH3D}
			\left\{
			\begin{split}
				&u(x,y,z;t=0)=\max\bigg(tanh\frac{0.26-R_1}{2\varepsilon},\tanh\frac{0.26-R_2}{2\varepsilon}\bigg),\\
				&R_1=\sqrt{(x-0.5)^2+(y-0.5)^2+(z-0.66)^2},\\
				&R_2=\sqrt{(x-0.5)^2+(y-0.5)^2+(z-0.34)^2}.\\
			\end{split}
			\right.
		\end{equation}\par
		In this experiment, $10,000$ initial points, $1,600$ boundary points, and $20,000$ residual points are used for training. Finally, we set the mass loss weight to be $10$ for the first subinterval and $1$ for the others. The weight $\lambda_u$ of the initial condition is set to $100$, and all others are $1$. To train the neural network, we first use the Adam optimizer with a learning rate of $0.0005$, then the L-BFGS-B optimizer to fine-tune the neural network. \par
		Adaptive resampling in space was performed every 5000 steps of training, and each time 1,000 residual points with large residuals of the equation were chosen and added to the initial set of residual points. This procedure is performed three times for each sub-time interval. The mapping proposed in Section \ref{FH_2D} is added for the output of the network to prevent the gradient explosion due to the singularity. The time domain is divided into four segments, namely with lengths of 0.02, 0.18, 0.4, and 0.4.\par
\begin{figure}[H]
			\centering
			\includegraphics[width=1\textwidth]{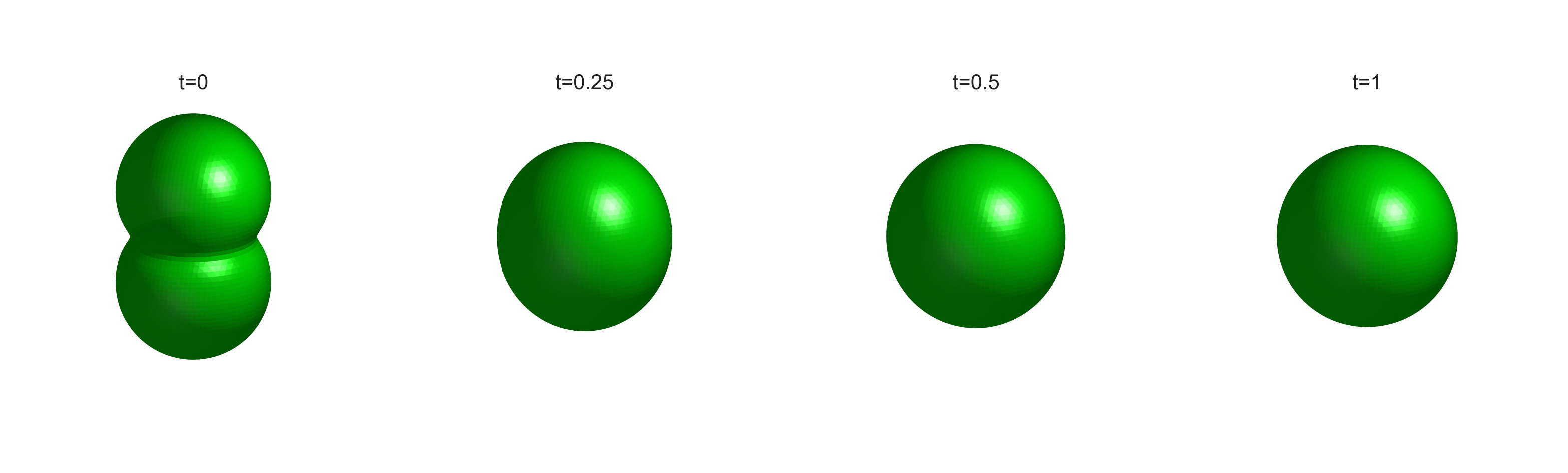}
			\caption{Predicted solutions at $t=0,0.25,0.5,1$. The lengths of the time subinterval are $0.02, 0.18, 0.4$ and $0.4$, respectively. The 3D surface is plotted by points in 3D space with the value $0$.}
			\label{FH3D_result}
		\end{figure}
\begin{figure}[H]
			\centering
			\subfloat[Energy]   
			{
				\includegraphics[width=0.45\textwidth]{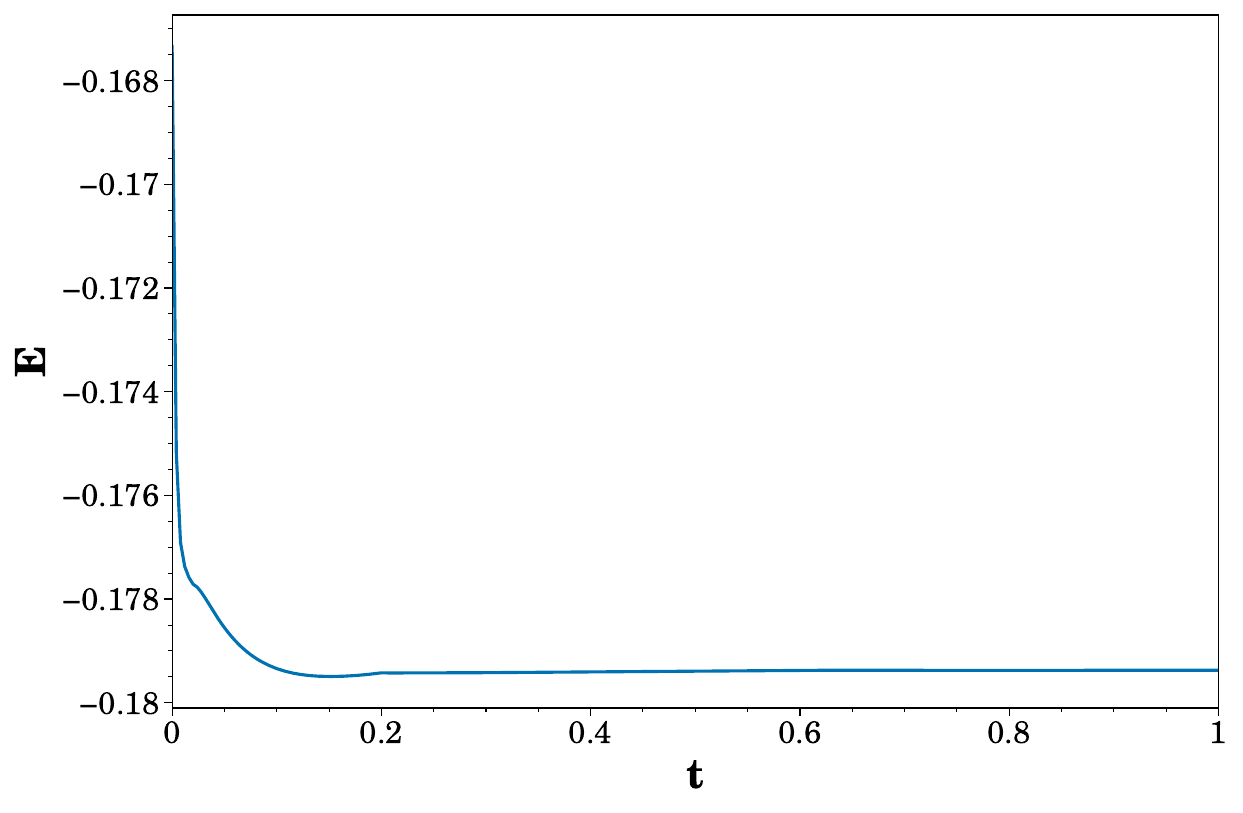}
				\label{FH3D_energy}
			}
			\subfloat[Mass]
			{
				\includegraphics[width=0.45\textwidth]{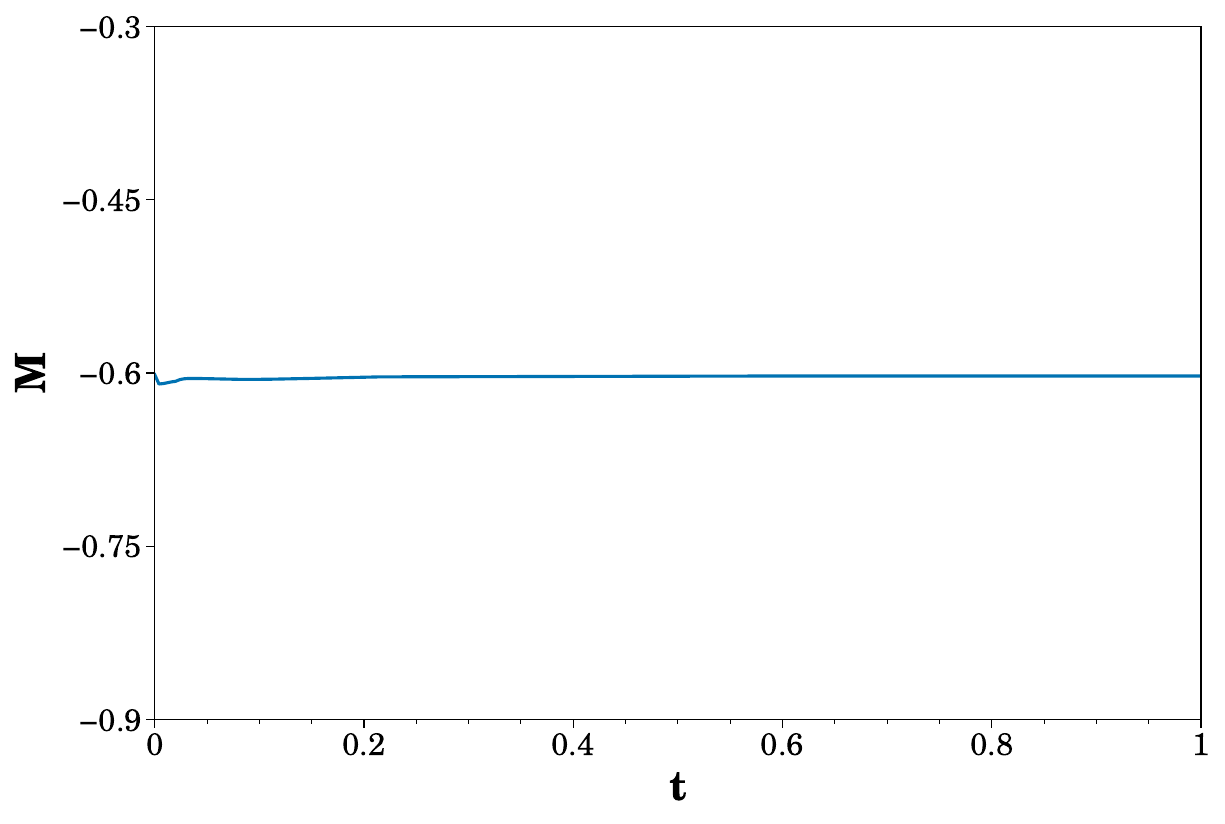}
				\label{FH3D_mass}
			}
			\caption{The numerical energy(E) and mass(M) of the predicted solutions $\hat{u}_\theta(\bm{x},t)$.}
			\label{FH3D_EM}
		\end{figure}

		The predicted solutions at different times ($t=0,0.25,0.5,1$) are shown in Fig.~\ref{FH3D_result}, showing the merging process from two balls into one. The reference solutions are computed by using a spatial discretization mesh of $64\times 64\times 64$ and uniform time subinterval of length $\Delta t=1e^{-5}$. The magnitude of the absolute errors is similar to one of the two-dimensional case, namely O(1).  Based on the predicted solutions above, the numerical energy and mass within $t\in[0,1]$ are plotted in Fig.~\ref{FH3D_EM}\subref{FH3D_energy} and Fig.~\ref{FH3D_EM}\subref{FH3D_mass}. The maximum mass error is 6.85e-03, which is relatively reasonable.
	
        \vspace{-0.3cm}
		The energy of the system first decreases and then tends to be stable with time. The predicted solutions preserve mass conservation. Similar results are obtained from the simulation of the Cahn-Hilliard equation with Ginzburg-Landau potential in three dimensions. From the above simulation, we conclude that the algorithm achieves a satisfactory level of simulation for the 3D Cahn-Hilliard equation.
		\subsubsection{System of Cahn-Hilliard equations in two dimensions}
		Next, we simulate the two-dimensional  system of Cahn-Hilliard equations with coupling terms, which is used to describe the confined binary blended polymers\cite{Hirai2019,Ji2022,Avalos2016}. The energy functional of the system is the following expression:
		\begin{equation}\label{E_CHS}
			E(u,v)=\int_{\Omega}\bigg(\frac{\varepsilon_u^2}{2}|\nabla u|^2+\frac{\varepsilon_v^2}{2}|\nabla v|^2+W(u,v)\bigg)d\bm{x},
		\end{equation}
		where
		\begin{equation}\label{W_CHS}
			W(u,v)=\frac{(u^2-1)^2}{4}+\frac{(v^2-1)^2}{4}+\alpha uv+\beta uv^2.
		\end{equation}
		It is a coupled system of two Cahn-Hilliard equations, where $u$ and $v$ are the order parameters distinguishing the polymer-rich and solvent-rich regions and the two polymer components, respectively. Parameters $\varepsilon_u$ and $\varepsilon_v$ control the width of the interface between macrophase and microphase. In equation (\ref{W_CHS}), the coupling parameter $\alpha$ controls the hydrophilicity of the microscopic phase, $\beta$ controls the formation of the restricted interface. For the energy functional(\ref{E_CHS}), two coupled Cahn-Hilliard equations are obtained under the $H^{-1}$ gradient flow,
		as follows:
		\begin{equation}\label{FH_example_2D}
			\left\{
			\begin{split}
				&u_t=\Delta(-\varepsilon_u^2\Delta u+u^3-u+\alpha v+\beta v^2),\\
				&v_t=\Delta(-\varepsilon_v^2\Delta v+v^3-v+\alpha u+2\beta uv),\\
			\end{split}
			\right.
		\end{equation}
		with periodic boundary conditions. The solution domain is $\Omega:=[-1,1]^2$, $t\in[0,1]$, and $\varepsilon_u=\varepsilon_v=0.05$. The coupling parameters are set as $\alpha=0$, $\beta=-0.5$ , meaning that the two components of the microphase have the same hydrophilicity. For this case, we choose the initial condition of $u$ as a circle and the initial condition of $v$ as two circles along the diagonal of the solution domain:
		\begin{equation}\label{initial_condition_CHS1}
			\left\{
			\begin{split}
				&u(x,y;t=0)=\tanh\bigg(\frac{0.6-\sqrt{x^2+y^2}}{2\varepsilon_u}\bigg),\\
				&v(x,y;t=0)=\tanh\bigg(\frac{0.4-R_1}{2\varepsilon_v}\bigg)-\tanh\bigg(\frac{0.4-R_2}{2\varepsilon_v}\bigg),\\
				&R_1=\sqrt{(x+0.25)^2+(y+0.25)^2},\\
				&R_2=\sqrt{(x-0.25)^2+(y-0.25)^2}.
			\end{split}
			\right.
		\end{equation}\par
		In this example, we use a neural network consisting of 10 hidden layers, with 50 neurons per layer. The above system is transformed into a system of four second order differential equations. Then we set the relevant residual loss functions. The number of sampling points for $u$ and $v$ is the same, with 10000 initial points, 1600 boundary points and 20000 residual points. The lengths of the time subintervals are $0.05, 0.15, 0.4, 0.2$ and $0.2$, and other settings are the same as before.\par
\begin{figure}[H]
			\centering         
			\subfloat[Predicted solutions $\hat{u}_\theta$ at $t=0,0.25,0.5,1$]   
			{
				\includegraphics[width=1\textwidth]{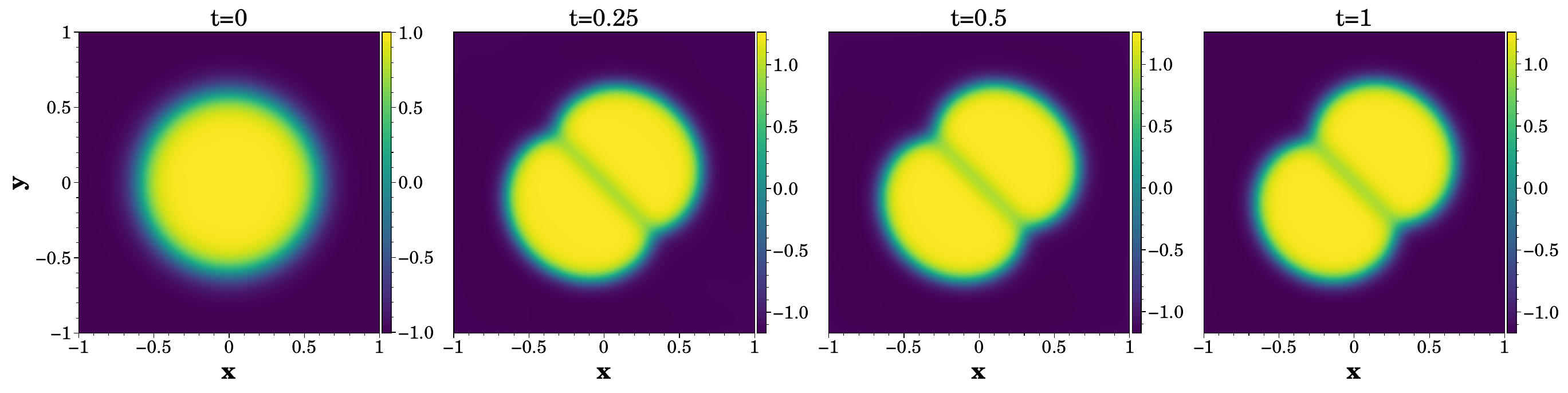}
				\label{result_chs_u1}
			}\\
			\subfloat[Predicted solutions $\hat{v}_\theta$ at $t=0,0.25,0.5,1$]
			{
				\includegraphics[width=1\textwidth]{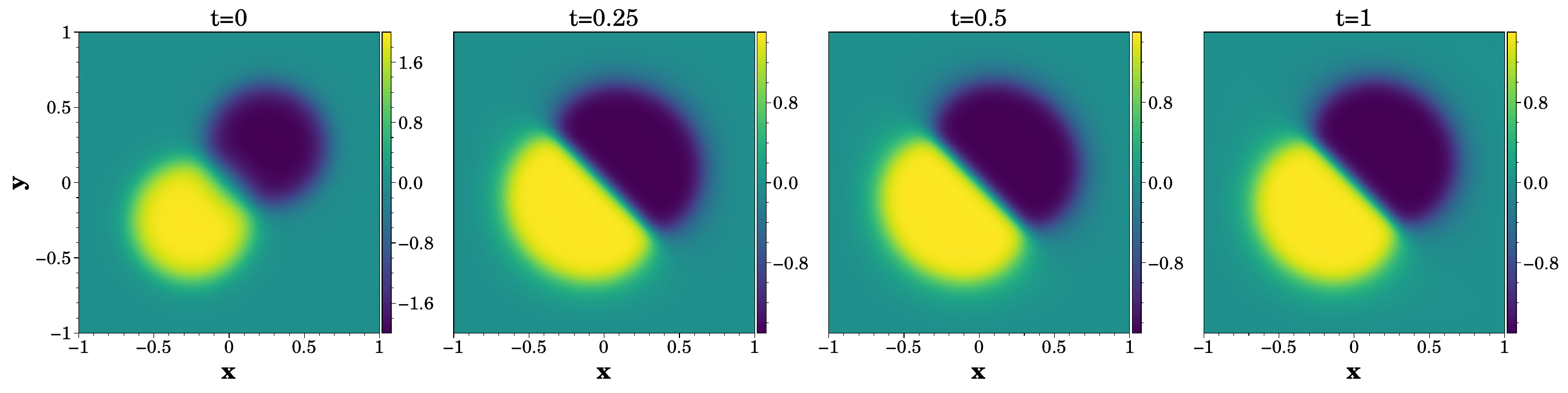}
				\label{result_chs_v1}
			}
			\caption{Predicted solutions $\hat{u}_\theta$ and $\hat{v}_\theta$ obtained by adaptive PINNs with time subintervals of lengths $0.05, 0.15, 0.4, 0.2$ and $0.2$.}
			\label{CHS_double}
		\end{figure}
\begin{figure}[H]
			\centering         
			\subfloat[Predicted solutions $\hat{u}_\theta$ at $t=0,0.25,0.5,1$]   
			{
				\includegraphics[width=1\textwidth]{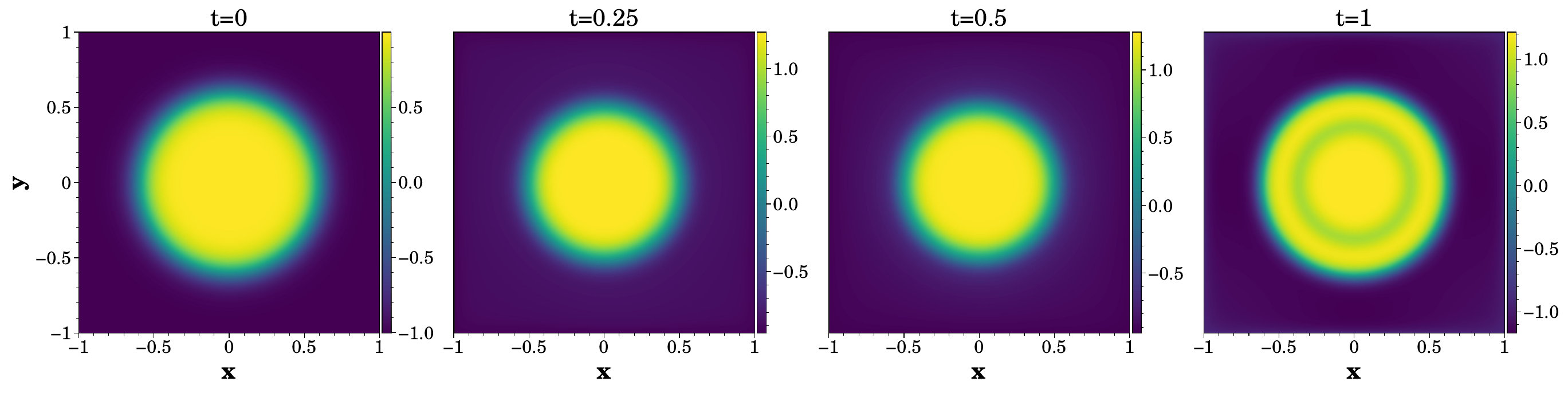}
				\label{result_chs_u2}
			}\\
			\subfloat[Predicted solutions $\hat{v}_\theta$ at $t=0,0.25,0.5,1$]
			{
				\includegraphics[width=1\textwidth]{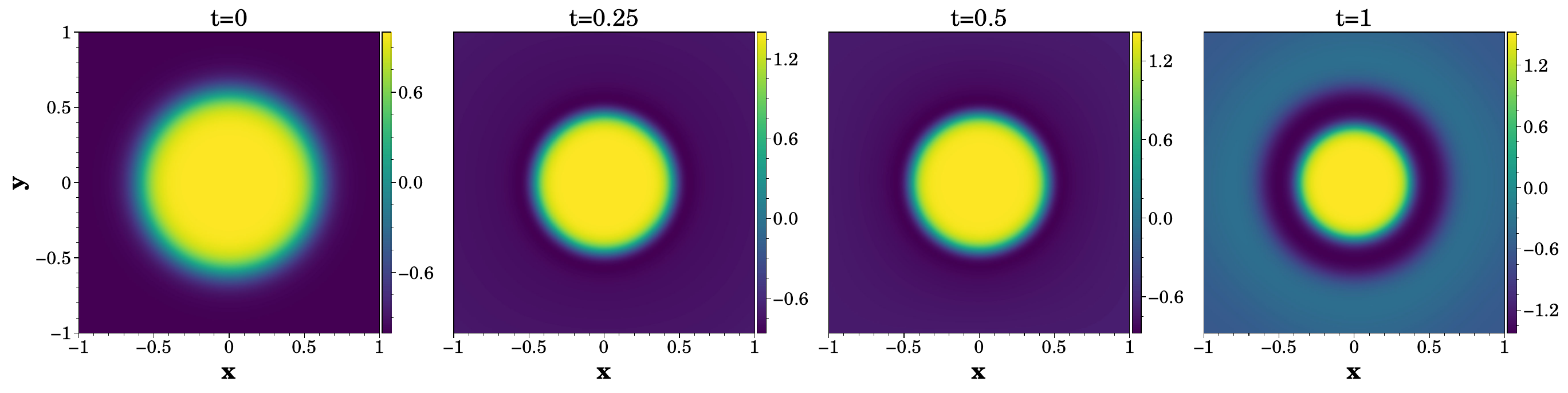}
				\label{result_chs_v2}
			}
			\caption{Predicted solutions $\hat{u}_\theta$ and $\hat{v}_\theta$ obtained by adaptive PINN with time subintervals of lengths $0.01, 0.04, 0.15, 0.4$ and $0.4$.
			}
			\label{CHS_onion}
		\end{figure}

		Fig.~\ref{CHS_double}\subref{result_chs_u1} shows the transition of the confined interface between polymers and solution. The structure of polymers, represented by the yellow region, changes from a circle to two semicircles. In Fig.~\ref{CHS_double}\subref{result_chs_v1}, the yellow and blue regions in the centre have values of $1$ and $-1$, representing the two kinds of polymers. The external region has a value of $0$, representing the solution. Eventually, each polymer occupies one-half of the confined domain, forming a half-sphere structure in the interior. 
		Next, we simulate the case of hydrophilicity, namely $\alpha=0.05$ and $\beta=-0.5$. Other parameters remain the same as above. In this case, the component represented by $v=-1$ has a stronger hydrophilicity. We choose circular initial conditions for both $u$ and $v$, as follows:
		\begin{equation}\label{initial_condition_CHS2}
			\left\{
			\begin{split}
				&u(x,y;t=0)=\tanh\bigg(\frac{0.6-\sqrt{x^2+y^2}}{2\varepsilon_u}\bigg),\\
				&v(x,y;t=0)=\tanh\bigg(\frac{0.6-\sqrt{x^2+y^2}}{2\varepsilon_v}\bigg).
			\end{split}
			\right.
		\end{equation}
		The lengths of time subintervals are set as $0.01, 0.04, 0.15$, $0.4$ and $0.4$, with other hyperparameters kept the same as the above experiment. The predicted solutions $\hat{u}_\theta$ and  $\hat{v}_\theta$ at different time are shown in Fig.~\ref{CHS_onion}\subref{result_chs_u2} and Fig.~\ref{CHS_onion}\subref{result_chs_v2}. \par
		Fig.~\ref{CHS_onion}\subref{result_chs_u2} and Fig.~\ref{CHS_onion}\subref{result_chs_v2} show the phase separation of $u$ and $v$, respectively. The interior is divided into yellow and dark blue regions, where the dark blue region surrounds the yellow region, similar to the structure of an onion. This indicates that the component represented by the dark blue has a stronger hydrophilicity.\par
		All the above results satisfy the property of mass conservation and the maximum absolute errors are of $O(1)$ and reasonable. From the simulation of these two complex cases, it can be seen that the mass-preserving spatio-temporal adaptive PINN that we have developed can well simulate high-dimensional or complex systems of equations.
        
\subsection{Comparison with other PINN variants}
To further demonstrate the effectiveness of our method, we compare it against two representative baselines: conservative PINN (cPINN) \cite{Jagtap2020cPINN} and PINN with adaptive activation functions (APINN) \cite{Jagtap2020Adaptive}. The comparison is performed on the 2D Cahn-Hilliard benchmark with the Ginzburg-Landau potential.

The results are summarized in Table~\ref{tab:comparison}. From this table, our method (STAPINN) achieves a lower maximum mass errors and a lower relative $L^2$ errors at $t=0.25$ compared to both baselines. We also discuss the trade-offs in model design and scalability. We note that XPINN \cite{Jagtap2020XPINN} and adaptive domain-decomposition approaches address different aspects of PINN improvement; a detailed comparison with these methods is left as future work.

    \begin{table}[H]
	\caption{Cahn-Hilliard equation with Ginzburg-Landau potential: relative $L^2$ errors and maximum mass error of three methods}
		\label{tab:comparison}
		\centering
		\def\temptablewidth{0.8\textwidth}
			\begin{tabular*}{\temptablewidth}{@{\extracolsep{\fill}}llcccc}
				\hline
				&Time&STAPINN&cPINN&APINN&\\ 
				\hline
                &0&1.39e-03&7.07e-04&3.15e-04&\\
				&0.25&8.36e-03&1.16e-01&6.28e-02&\\
				\hline
				&mass error&1.01e-03&5.86e-02&1.15e-01&\\
				\hline
			\end{tabular*}\\ 
	\end{table}

		\section{Conclusion}
		In this paper, we have introduced the mass-preserving spatio-temporal adaptive PINN that significantly improves the accuracy of PINN for solving high-order phase field equations with strong nonlinearity, singularities, and mass conservation property. Specifically, we add a mass constraint to the loss function and combine it with adaptive sampling methods in space and time. 
        This paper provides a detailed description of the process of adjusting the lengths of time subintervals based on the rate of energy decrease and modifying their weights according to the mass error; spatially, residual-based resampling focuses collocation points near moving interfaces.

		We conduct several numerical experiments on the Cahn-Hilliard equations with Ginzburg-Landau potential and Flory-Huggins potential. Numerical results show that the proposed modified PINN produces more accurate approximate solutions compared to the baseline PINN with uniform time subintervals, while also ensuring mass conservation. 
        A sensitivity analysis on the energy threshold further confirms the robustness of our time partitioning strategy.
        Comparisons with other representative PINN variants, including Conservative PINN and PINN with adaptive activation functions, further confirm the superiority of our method in terms of accuracy and mass preservation. Finally, our simulation results for the three-dimensional problem and the system of Cahn-Hilliard equations further indicate the effectiveness and accuracy of the approach. In summary, our modified PINN method is capable of solving high-order phase field equations with nonlinearity and singularity.\par

In future work, we plan to extend our framework to other complex phase-field and gradient-flow models, develop a more theoretically guided energy-partition strategy, and explore automatic hyperparameter tuning for adaptive sampling and mass conservation enforcement. More importantly, we aim to move beyond soft constraints by designing neural network structures that inherently preserve physical properties, such as mass conservation and solution boundedness, without adding extra penalty terms \cite{chen2025}. 
Exploring alternative network architectures, such as tensor neural networks (TNNs), is also a promising direction to further improve accuracy.
We also note that GNN interpretability methods, such as the functional-semantic mapping \cite{Raj2024}, may provide useful insights for analyzing the learned partitions in our adaptive PINN framework. Thus, we will work toward a more accurate and efficient deep learning method for solving complex phase-field equations with guaranteed conservation laws and other physical structures.

\section*{Acknowledgements}
Q. Huang was partially supported by the National Natural Science Foundation of China (No. 12371385). Z. Xu was partially supported by the National Natural Science Foundation of China (No. 12201021).
		
		
		\bibliography{reference_short}
		\bibliographystyle{elsarticle-num} 
		
	\end{document}